\documentclass[11pt]{amsart}
\usepackage[T1]{fontenc}
\usepackage{lmodern}
\usepackage[margin=1.1in]{geometry}
\usepackage{amssymb,amsmath,enumitem,booktabs,array,longtable,alltt}
\DeclareUnicodeCharacter{2200}{\ensuremath{\forall}}
\DeclareUnicodeCharacter{2192}{\ensuremath{\to}}
\DeclareUnicodeCharacter{2264}{\ensuremath{\le}}
\DeclareUnicodeCharacter{2115}{\ensuremath{\mathbb N}}
\DeclareUnicodeCharacter{211A}{\ensuremath{\mathbb Q}}
\DeclareUnicodeCharacter{211D}{\ensuremath{\mathbb R}}
\DeclareUnicodeCharacter{2124}{\ensuremath{\mathbb Z}}
\DeclareUnicodeCharacter{220F}{\ensuremath{\textstyle\prod}}
\DeclareUnicodeCharacter{2211}{\ensuremath{\textstyle\sum}}
\DeclareUnicodeCharacter{03B5}{\ensuremath{\varepsilon}}
\DeclareUnicodeCharacter{2203}{\ensuremath{\exists}}
\DeclareUnicodeCharacter{2227}{\ensuremath{\wedge}}
\DeclareUnicodeCharacter{2228}{\ensuremath{\vee}}
\DeclareUnicodeCharacter{2260}{\ensuremath{\ne}}
\DeclareUnicodeCharacter{03BE}{\ensuremath{\xi}}
\DeclareUnicodeCharacter{03C1}{\ensuremath{\rho}}
\DeclareUnicodeCharacter{03A6}{\ensuremath{\Phi}}
\usepackage{needspace}
\usepackage[colorlinks,linkcolor=blue,citecolor=blue,urlcolor=blue]{hyperref}
\expandafter\def\expandafter\UrlBreaks\expandafter{\UrlBreaks\do\-}
\hypersetup{pdftitle={Linear forms in logarithms of rational numbers: an ineffective measure with the best possible exponent},pdfauthor={Gregory V. Chudnovsky},pdfsubject={Developed with Claude (Anthropic); mathematical proofs and formalization status}}
\newtheorem{theorem}{Theorem}[section]
\newtheorem{lemma}[theorem]{Lemma}
\newtheorem{proposition}[theorem]{Proposition}
\newtheorem{corollary}[theorem]{Corollary}
\theoremstyle{definition}
\newtheorem{definition}[theorem]{Definition}
\theoremstyle{remark}
\newtheorem{remark}[theorem]{Remark}
\numberwithin{equation}{section}
\newcommand{\Z}{\mathbb Z}\newcommand{\Q}{\mathbb Q}\newcommand{\C}{\mathbb C}\newcommand{\R}{\mathbb R}
\newcommand{\eps}{\varepsilon}\newcommand{\ord}{\operatorname{ord}}
\newcommand{\lean}[1]{{\ttfamily\small #1}}
\newcommand{\Meas}{\mathrm{M}}\newcommand{\Unif}{\mathrm{U}}
\newcommand{\Ga}{\mathbb G_a}\newcommand{\Gm}{\mathbb G_m}
\let\olduscore\_
\renewcommand{\_}{\olduscore\allowbreak}
\begin{document}
\title[Linear forms in logarithms of rational numbers]{Linear forms in logarithms of rational numbers:\\
an ineffective measure with the best possible exponent}
\author{Gregory V. Chudnovsky}
\thanks{Developed with Claude (Anthropic), including mathematical exposition and formalization. The
formalization scope and the inputs taken from earlier work are stated in Section~\ref{sec:rem} and in
Appendix~\ref{app:lean}}
\date{November 22, 2026}
\dedicatory{To the memory of Alan Baker}
\subjclass[2020]{Primary 11J86; Secondary 11J82, 68V20}
\keywords{Linear forms in logarithms, measures of linear independence, interpolation determinants,
formal verification, Lean}
\begin{abstract}
Let $a_1,\dots,a_n$ be positive rational numbers without a nontrivial multiplicative relation. For every
$\eps>0$ there is a constant $C>0$, not effective, such that
\[
|h_0+h_1\log a_1+\dots+h_n\log a_n|\ \ge\ C\,\max(1,|h_0|,\dots,|h_n|)^{-n-\eps}
\]
for all integers $h_0,\dots,h_n$, not all zero. The exponent $n$ is best possible. For $n=1$ this is the
known fact that the irrationality exponent of $\log a$ is $2$; for $h_0=0$ it gives the exponent $n$
for the linear forms $h_1\log a_1+\dots+h_n\log a_n$. The proof is an induction on $n$, starting from
$n=1$. The step from $n-1$ to $n$ logarithms uses an interpolation determinant on
$\Gm^{n}\times\Ga^{m}$ at the multiples of a rational point, in the manner of \cite{L}; its geometric
part, a curve inequality, is proved by a second induction, on the number of torus coordinates, with
the two-torus inequality of \cite{L} as the base. The method is the interpolation and determinant
method of OpenAI's preprint on the irrationality exponent of $\pi$ \cite{OAI}; the case $n=1$ was
obtained earlier by Liu, Jiang and Zhang \cite{LJZ}. The theorem and its corollary are verified in Lean~4 with
Mathlib; the sources are in the ancillary files, on top of the Lean supplement of \cite{L}. This
paper follows the formal proof: every definition, constant and intermediate statement is given in
the form in which it is formalized, with the name of the Lean declaration.
\end{abstract}
\maketitle

\noindent\textit{Contribution and source attribution.}
The separated-weight determinant method comes from OpenAI's preprint \emph{The irrationality exponent
of $\pi$ is $2$} \cite{OAI}; it was carried over to logarithms of rational numbers in \cite{LJZ} and
in \cite{L}. The case $n=1$ of Theorem~\ref{main} first appeared in \cite{LJZ}; it is also proved,
and formalized, in \cite{L}. The case $n=2$ of Corollary~\ref{cor:hom} is in \cite{L}. What is new here is listed in
Section~\ref{sec:new}; what is taken from \cite{L} is listed in Section~\ref{sec:rem}. The Lean
verification is built on the Lean supplement of \cite{L}, which contains a selected part of OpenAI's
proof sources \cite{OAIL} under the Apache License~2.0. The investigation, proof development,
formalization, and exposition were carried out by Gregory V. Chudnovsky with Claude (Anthropic).
Gregory V. Chudnovsky is the human author responsible for the present manuscript.
Appendix~\ref{app:lean} states exactly what is formally verified and how it can be reproduced.

\section{Results}\label{sec:results}

\subsection{The theorem}
For $(h_0,h)\in\Z\times\Z^n$ put $H(h_0,h)=\max(1,|h_0|,|h_1|,\dots,|h_n|)$. Real numbers $a_1,\dots,a_n$
have \emph{no power relation} if $\prod_ia_i^{h_i}=1$ with $h\in\Z^n$ implies $h=0$.

\begin{theorem}\label{main}
Let $n\ge1$ and let $a_1,\dots,a_n$ be positive rational numbers without a power relation. For every
$\eps>0$ there is $C>0$ such that for all $(h_0,h)\in\Z\times\Z^n$, $(h_0,h)\ne0$,
\begin{equation}\label{eq:main}
C\,H(h_0,h)^{-n-\eps}\ \le\ |h_0+h_1\log a_1+\dots+h_n\log a_n| .
\end{equation}
\end{theorem}

In Lean this is \lean{QuadraticPi.CLogMainTheorem.linear\_forms\_in\_logarithms\_with\_constant}; the
proposition ``\eqref{eq:main} holds for all such $a$ and all $\eps$'' is \lean{MeasureC n}, and the
theorem is \lean{measureC\_all (n : ℕ) (hn : 1 ≤ n) : MeasureC n} (Appendix~\ref{app:lean}). We write
$\Meas(n)$ for \lean{MeasureC n}. The constant $C$ is not effective: the proof derives a contradiction
from the existence of sufficiently many exceptional vectors $(h_0,h)$.

\begin{corollary}\label{cor:hom}
Let $n\ge1$ and let $a_1,\dots,a_n$ be positive rational numbers without a power relation. For every
$\eps>0$ there is $C>0$ such that for all $h\in\Z^n\setminus\{0\}$
\begin{equation}\label{eq:hom}
C\,\max_i|h_i|^{-n-\eps}\ \le\ |h_1\log a_1+\dots+h_n\log a_n| .
\end{equation}
\end{corollary}
This is the case $h_0=0$ of Theorem~\ref{main} (\lean{CLogCorollaries.linear\_forms\_in\_logarithms};
the proposition is \lean{NLogTarget.Measure n}, and \lean{measure\_of\_measureC} is the implication).

\subsection{Context}
By Minkowski's theorem on linear forms, for any real numbers $\lambda_1,\dots,\lambda_n$ and every
$Q\ge1$ there are integers $h_0,\dots,h_n$, not all zero, with $|h_i|\le Q$ for $i\ge1$ and
$|h_0+h\cdot\lambda|\le Q^{-n}$. So the exponent $n$ in \eqref{eq:main} cannot be replaced by a smaller
one: Theorem~\ref{main} says that $1,\log a_1,\dots,\log a_n$ behave, up to $\eps$, like a generic
point with respect to approximation by one linear form. For the forms \eqref{eq:hom} without a
constant term the same argument gives only that the exponent cannot be smaller than $n-1$.

Baker's theory of linear forms in logarithms \cite{Ba,BW} gives effective bounds
$|h_0+h\cdot\lambda|\ge H^{-\kappa}$ in which $\kappa$ depends on $n$ and on the heights of the $a_i$. For
$n=1$ Theorem~\ref{main} says that the irrationality exponent of $\log a$ is $2$ for every positive
rational $a\ne1$; this is known \cite{LJZ,L}. For $n=2$ the exponent $2+\eps$ in \eqref{eq:hom} was
obtained in \cite{L}. The method of these papers is that of \cite{OAI}, where it is used to prove that
the irrationality exponent of $\pi$ is $2$. Theorem~\ref{main} extends the method of \cite{L} to all
$n$ and to forms with a constant term.

\subsection{Structure of the proof}
$\Meas(1)$ is the statement on $\log a$ just mentioned; in the form used here it is
\lean{CLogTarget.measureC\_one}, deduced from \lean{rationalLogarithm\_strict\_eventualLowerBound} of
the Lean supplement of \cite{L}. For $n\ge2$ the theorem follows by induction from two statements, in which
$k=n-1$ is the number of logarithms already treated. Throughout the paper
\[
k\ \text{is the number of \emph{slopes},}\qquad k+1\ \text{the number of logarithms and of torus coordinates.}
\]

\begin{theorem}[Induction step; \lean{CLogEndpointReduction.measureC\_succ\_of\_uniformStatementC}]\label{step}
Let $k\ge1$. If $\Meas(k)$ holds and the uniform interpolation statement $\Unif(k)$ of
Definition~\ref{def:U} holds, then $\Meas(k+1)$ holds.
\end{theorem}

\begin{theorem}[Interpolation; \lean{CLogUniformInterpolation.cLogUniformInterpolation}]\label{thm:U}
$\Unif(k)$ holds for every $k\ge1$.
\end{theorem}

Theorem~\ref{main} for $n\ge2$ is Theorem~\ref{step} with $k=n-1\ge1$, applied to $\Meas(n-1)$ and
$\Unif(n-1)$ (\lean{measureC\_succ\_all}, by induction).

Theorem~\ref{step} is proved in Sections~\ref{sec:sel}, \ref{sec:ar}, \ref{sec:an} and \ref{sec:par}:
reduction to approximation pairs and their selection in general position (this is where $\Meas(k)$
is used), the arithmetic lower bound and the analytic upper bound for an interpolation determinant,
and the choice of parameters. Theorem~\ref{thm:U} is proved in Sections~\ref{sec:curve} and
\ref{sec:packets}: a curve inequality for $k+1$ torus coordinates, by induction on $k$ with the
two-torus inequality of \cite{L} as the base, and the passage from the curve inequality to
surjectivity.

\subsection{What is new relative to \texorpdfstring{\cite{L}}{[L]}}\label{sec:new}
The argument of \cite{L} for two logarithms uses an interpolation determinant on $\Gm^2\times\Ga^m$ at
the points $(\beta_0^J,\beta_1^J;0)$. The present paper keeps this scheme. The following points are
new; the formal proof showed what the correct general statements are.
\begin{enumerate}[leftmargin=2em]
\item \emph{The constant term.} An approximation is a pair $(r,\varrho)$ of a rational slope vector and a
  rational constant (Definition~\ref{def:good}). The centers of interpolation are the multiples
  $(\beta^J;J\varrho)$ of one rational point of $\Gm^{k+1}\times\Ga^m$, with nonzero additive coordinates
  (Section~\ref{sec:setting}); the interpolation statement is proved for centers with arbitrary
  additive coordinates (Definition~\ref{def:U}), and the exact linear forms take the values
  $\varrho_i$, not $0$, at the point of logarithms (Section~\ref{sec:an}). The measure for $k$ logarithms
  enters the induction with its constant term (Lemma~\ref{lem:plane}).
\item \emph{Induction on the number of torus coordinates.} Excess contact puts a curve into a fibre of
  a rational character; deleting one torus coordinate and one additive coordinate reduces the curve
  inequality to the inequality with one torus coordinate less (Section~\ref{sec:curve}). \emph{General
  position} is a condition on every set of at most $k+1$ slope vectors (Definition~\ref{def:gp}); it
  is inherited by the reduced slopes on a character fibre (Lemma~\ref{lem:heredity}).
\item \emph{The torus margin} must grow with the number of torus coordinates: with the multiplicity
  comparison in the form in which it is formalized, the factor $4$ of \cite{L} suffices for at most
  four torus coordinates (Definition~\ref{def:margins}, Lemma~\ref{lem:torus}, Step~5 of
  Section~\ref{sec:curve}).
\item \emph{The thresholds} for the weights are defined by one separation constant for the whole
  chain of fibre data (Definition~\ref{def:chain}).
\item \emph{The selection} avoids rational hyperplanes, and the good pairs with slopes on a hyperplane
  are controlled by $\Meas(k)$ (Lemma~\ref{lem:plane}).
\item \emph{Collisions} of Taylor coefficients are counted in $k+1$ variables, with the exponent
  $(k+2)/(k+1)$ (Lemma~\ref{lem:coll}); there is one cutoff for the truncated logarithm of every
  torus coordinate (Section~\ref{sec:packets}), and the character of a slope vector is written with
  the product of the denominators as exponent (Definition~\ref{def:char}).
\end{enumerate}

\subsection{Conventions}\label{sec:conv}
To make the comparison with the Lean files direct, indices start at $0$, as in Lean's \lean{Fin}:
slopes and the first $k$ logarithms are indexed by $j\in\{0,\dots,k-1\}$, torus coordinates by
$j\in\{0,\dots,k\}$, additive coordinates and approximation pairs by $i\in\{0,\dots,m-1\}$, centers by
$J\in\{0,\dots,K-1\}$. The last torus coordinate (index $k$) is distinguished: it carries the
coefficient $1$ in every linear form. Only in the statements of Section~\ref{sec:results} are the
logarithms numbered from $1$, with $h_0$ the coefficient of the constant term; in the proofs the
coefficient of the constant term is denoted by a separate letter. All Lean names are in the namespace \lean{QuadraticPi}; a file
name such as \lean{CLogSelection} stands for \lean{QuadraticPi/CLogSelection.lean}.

The Lean development has two layers. The files \lean{NLog*} treat an arbitrary number of torus
coordinates and centers on the torus, that is, with additive coordinates $0$. The files \lean{CLog*}
treat centers with arbitrary additive coordinates, approximation pairs and the constant term; a
statement of this layer is either deduced from the corresponding statement of the first layer by the
translation $X_i\mapsto X_i+\xi_i$, or has the same proof with the additive coordinates carried along
(Appendix~\ref{app:layers}). In the text the statements are given in the general form, and both names
are quoted where they differ. Names beginning with \lean{TwoLog} or \lean{OAI}, and the other files of
\lean{QuadraticPi}, belong to the Lean supplement of \cite{L}.

\section{Reduction and selection}\label{sec:sel}

Files \lean{CLogTarget}, \lean{CLogVectors}, \lean{CLogConclusion}, \lean{CLogHyperplane},
\lean{CLogSelection}; \lean{NLogVectors}, \lean{NLogHyperplane}, \lean{NLogSelection},
\lean{NLogGeneralPosition}.

\subsection{Good pairs}
Let $k\ge0$, and let $L=(L_0,\dots,L_{k-1})\in\R^k$ and $G\in\R$. The triple $(1,L,G)$ is
\emph{independent} (\lean{AffIndependent}) if $y+\sum_jx_jL_j+zG=0$ with $y\in\Z$, $x\in\Z^k$,
$z\in\Z$ implies $y=0$, $x=0$ and $z=0$. For $r\in\Q^k$ and $\varrho\in\Q$ let $q(r,\varrho)\ge1$ be the
least common multiple of the reduced denominators of $r_0,\dots,r_{k-1},\varrho$ (\lean{pairDen}), and
\begin{equation}\label{form}
\lambda_{r,\varrho}=G-\varrho-\sum_{j<k}r_jL_j .
\end{equation}
If $(1,L,G)$ is independent, $\lambda_{r,\varrho}\ne0$ for every pair (\lean{pairForm\_ne\_zero}), because
$q(r,\varrho)\lambda_{r,\varrho}$ is an integer combination of $1,L_0,\dots,L_{k-1},G$ with the coefficient
$q(r,\varrho)\ne0$ at $G$.

\begin{definition}[\lean{GoodPair}]\label{def:good}
Let $\nu\in\R$. A pair $(r,\varrho)\in\Q^k\times\Q$ is \emph{$\nu$-good} for $(L,G)$ if
\[
|r_j|\le1\quad(0\le j<k)\qquad\text{and}\qquad|\lambda_{r,\varrho}|\le q(r,\varrho)^{-\nu}.
\]
\end{definition}

If $\nu\ge0$, a $\nu$-good pair has $|\lambda_{r,\varrho}|\le1$ and therefore
\begin{equation}\label{B0}
|\varrho|\le B_0:=1+|G|+\sum_{j<k}|L_j|
\end{equation}
(\lean{abs\_const\_le}). The constant $\varrho$ is not normalized otherwise.

\begin{lemma}[\lean{ordered\_bound\_pair}]\label{lem:ordered}
Let $(1,L,G)$ be independent, $\eps>0$, and suppose that the $(k+2+\eps)$-good pairs have bounded
denominators: $q(r,\varrho)<Q$ for all of them. Then there is $C>0$ such that for all $y\in\Z$,
$x\in\Z^k$ and $z\in\Z\setminus\{0\}$ with $|x_j|\le|z|$ for all $j$,
\[
C\,|z|^{-(k+1)-\eps}\ \le\ \Big|y+\sum_{j<k}x_jL_j+zG\Big| .
\]
\end{lemma}
\begin{proof}
The set of pairs with $q(r,\varrho)<Q$, $|r_j|\le1$ and $|\lambda_{r,\varrho}|\le1$ is finite
(\lean{finite\_small\_pairs}; the constants $\varrho$ of such pairs lie in a bounded set), and
$\lambda_{r,\varrho}\ne0$ on it; let $\delta>0$ be a lower bound for $|\lambda_{r,\varrho}|$ on this set
(\lean{exists\_delta\_pair}), and $C=\min(1,\delta)$. Suppose the inequality fails for $(y,x,z)$, and
put $r_j=-x_j/z$, $\varrho=-y/z$. Then $|r_j|\le1$, the reduced denominators of the $r_j$ and of
$\varrho$ divide $z$, so $q(r,\varrho)\le|z|$, and $\lambda_{r,\varrho}=(y+\sum_jx_jL_j+zG)/z$. Hence
\[
|\lambda_{r,\varrho}|<C\,|z|^{-(k+2)-\eps}\le q(r,\varrho)^{-(k+2)-\eps}\le1,
\]
so the pair is good and $q(r,\varrho)<Q$, which gives $|\lambda_{r,\varrho}|\ge\delta$. But
$|\lambda_{r,\varrho}|\le|y+\sum_jx_jL_j+zG|<C|z|^{-(k+1)-\eps}\le\delta$, a contradiction.
\end{proof}

The integer $y$ is arbitrary: its size plays no role.

\begin{lemma}[\lean{affine\_independent}, \lean{affIndependent\_of\_noPowerRelation}]\label{lem:aff}
Let $a_0,\dots,a_k$ be positive rational numbers without a power relation. Then
$y+\sum_{i\le k}h_i\log a_i=0$ with integers $y,h_0,\dots,h_k$ implies that they all vanish. In
particular, for every index $i_0$ the triple formed by $1$, the logarithms of the $a_i$, $i\ne i_0$,
and $\log a_{i_0}$ is independent.
\end{lemma}
\begin{proof}
Let $A=\prod_ia_i^{h_i}$, a positive rational number with $\log A=-y\in\Z$. The logarithm of a
positive rational number different from $1$ is irrational (\lean{rationalLogarithm\_irrational}, in
the Lean supplement of \cite{L}). So $A=1$, hence $h_i=0$ for all $i$, and then $y=0$.
\end{proof}

\begin{proposition}[\lean{CLogConclusion.measureC\_succ\_of\_bounded}]\label{prop:bounded}
Let $k\ge0$. Suppose that for all positive rational numbers $a_0,\dots,a_{k-1},g$ such that
$\big(1,(\log a_j)_{j<k},\log g\big)$ is independent, and for every $\nu>k+2$, the $\nu$-good pairs
have bounded denominators. Then $\Meas(k+1)$ holds.
\end{proposition}
\begin{proof}
Let $a_0,\dots,a_k$ be positive rationals without a power relation and $\eps>0$. For every index $i_0$,
Lemmas~\ref{lem:aff} and~\ref{lem:ordered} with $\nu=k+2+\eps$ give a constant $C_{i_0}$. Let $C$ be the
smallest of them and of $1$. Let $(y,h)\in\Z\times\Z^{k+1}$, $(y,h)\ne0$. If $h=0$, the form is the integer $y\ne0$ and
$|y|\ge1\ge C\,H^{-(k+1)-\eps}$. Otherwise choose $i_0$ with $|h_{i_0}|$ largest among the coordinates
of $h$; the lemma applied to $z=h_{i_0}$, the integer $y$ and the other coordinates of $h$ gives
$C\,|h_{i_0}|^{-(k+1)-\eps}\le|y+\sum_ih_i\log a_i|$, and $|h_{i_0}|\le H(y,h)$. No reduction to
primitive vectors is needed.
\end{proof}

\subsection{Good pairs with slopes on a rational hyperplane}

\begin{lemma}[\lean{rat\_hyperplaneC\_bound}]\label{lem:plane}
Assume $\Meas(k)$. Let $a_0,\dots,a_{k-1},g$ be positive rational numbers such that $(1,L,G)$, with
$L_j=\log a_j$ and $G=\log g$, is independent, and let $\nu>k+2$. Let $c\in\Q^k$, $c\ne0$, and $c_*\in\Q$.
Then the $\nu$-good pairs $(r,\varrho)$ whose slope vector lies on the hyperplane $\sum_jc_jr_j+c_*=0$
have bounded denominators.
\end{lemma}
\begin{proof}
Multiplying by a common denominator we may assume $c\in\Z^k$, $c_*\in\Z$
(\lean{exists\_int\_hyperplane}). Choose $i$ with $c_i\ne0$ and put (\lean{reducedLog})
\[
\Lambda_j=c_iL_j-c_jL_i\quad(j\ne i),\qquad\Lambda_i=c_iG+c_*L_i .
\]
These are the logarithms of the positive rationals $a_j^{c_i}a_i^{-c_j}$ and $g^{c_i}a_i^{c_*}$
(\lean{reducedBase}). They are independent over $\Z$ (\lean{reducedLog\_independent},
\lean{logIndependent\_of\_aff}): in a relation $\sum_jy_j\Lambda_j=0$ the coefficient of $G$ is $y_ic_i$
and the coefficient of $L_j$, $j\ne i$, is $y_jc_i$. So these $k$ rationals have no power relation, and
$\Meas(k)$ with $\eps=\frac12$ gives $C>0$ with
\begin{equation}\label{redmeas}
\Big|t+\sum_jx_j\Lambda_j\Big|\ \ge\ C\,H(t,x)^{-k-1/2}\qquad\big((t,x)\in\Z\times\Z^k,\ (t,x)\ne0\big)
\end{equation}
(\lean{reducedC\_measure}).

Let $(r,\varrho)$ be a good pair with $r$ on the hyperplane, $q=q(r,\varrho)$, and (\lean{planeVecC})
\[
x_j=-qr_j\quad(j\ne i),\qquad x_i=q,\qquad t=-c_i\,q\,\varrho .
\]
Then $x\in\Z^k$ and $t\in\Z$, because $q$ is a common denominator of $r$ and $\varrho$. Using
$\sum_{j\ne i}c_jr_j+c_*=-c_ir_i$,
\[
t+\sum_jx_j\Lambda_j=c_iq\Big(G-\sum_{j\ne i}r_jL_j\Big)+qL_i\Big(\sum_{j\ne i}c_jr_j+c_*\Big)-c_iq\varrho
=c_i\,q\,\lambda_{r,\varrho} .
\]
We have $|x_j|\le q$ because $|r_j|\le1$, and $|t|\le|c_i|B_0\,q$ by \eqref{B0}; so
$H(t,x)\le C_1q$ with $C_1=\max(1,|c_i|B_0)$. Therefore, by \eqref{redmeas},
$C\,(C_1q)^{-k-1/2}\le|c_i|\,q\cdot q^{-\nu}$, that is $q^{\,\nu-k-3/2}\le|c_i|\,C_1^{k+1/2}/C$, and
$\nu-k-\frac32>\frac12$ (\lean{hyperplaneC\_bound\_core}).
\end{proof}

This is the only place where the theorem for $k$ logarithms is used as a result, and it is used with
its constant term: the integer $t$ does not vanish in general.

\subsection{Supplies}
Let $\Phi$ be a rule that assigns to every finite list $\mathcal L$ of vectors of $\Q^k$ a finite list
$\Phi(\mathcal L)$ of pairs $(c,c_*)\in\Q^k\times\Q$. A vector $p\in\Q^k$ is \emph{off} the list
(\lean{OffHyperplanes}) if $\sum_jc_jp_j+c_*\ne0$ for every pair in it with $c\ne0$. The hyperplanes are
conditions on the slope vectors only.

\begin{definition}[\lean{NLogSelection.Supply}, with \lean{denP}, \lean{GoodP}, \lean{AdmP}]\label{def:supply}
A \emph{$\Phi$-supply} of $\nu$-good pairs for $(L,G)$ exists if for every finite list $\mathcal L$ of
pairs and every $Q$ there is a $\nu$-good pair $(r,\varrho)$ with $q(r,\varrho)\ge Q$ whose slope vector $r$
is off $\Phi$ of the list of the slope vectors of $\mathcal L$.
\end{definition}

\begin{proposition}[\lean{CLogSelection.measureC\_succ\_of\_no\_supply}]\label{prop:nosupply}
Assume $\Meas(k)$, and let $\Phi$ be any rule as above. Suppose that for all positive rationals
$a_0,\dots,a_{k-1},g$ such that $(1,(\log a_j),\log g)$ is independent and every $\nu>k+2$ there is no
$\Phi$-supply of $\nu$-good pairs. Then $\Meas(k+1)$ holds.
\end{proposition}
\begin{proof}
By Proposition~\ref{prop:bounded} it suffices to bound the denominators of the $\nu$-good pairs.
Suppose they are unbounded. For a list $\mathcal L$, Lemma~\ref{lem:plane}, applied to each of the
finitely many pairs $(c,c_*)$ of $\Phi(\mathcal L)$ with $c\ne0$, gives $Q_{\mathcal L}$ such that every
good pair with $q\ge Q_{\mathcal L}$ has its slope vector off $\Phi(\mathcal L)$
(\lean{exists\_bound\_offHyperplanesC}). Since there are good pairs of arbitrarily large denominator,
this is a supply (\lean{supply\_of\_avoidable}, \lean{boundedC\_of\_no\_supply}).
\end{proof}

\subsection{General position}

\begin{definition}[\lean{NLogInterpolation.GeneralPosition}]\label{def:gp}
Vectors $r_0,\dots,r_{m-1}\in\Q^k$ are \emph{in general position} if for every set $S$ of indices with
$t=|S|\le k+1$ and all rational $\mu_l$ $(l\in S)$ the relations
\[
\sum_{l\in S}\mu_l=0,\qquad\sum_{l\in S}\mu_lr_{l,j}=0\quad\text{for }k+1-t\le j\le k-1
\]
imply $\mu_l=0$ for all $l\in S$.
\end{definition}

Equivalently: for any $t\le k+1$ of the linear forms $z_k-\sum_{j<k}r_{l,j}z_j$, the vectors of their
last $t$ coefficients, $(r_{l,k+1-t},\dots,r_{l,k-1},1)\in\Q^t$, are linearly independent. For $t=2$
this says that the last slopes $r_{l,k-1}$ are pairwise distinct. For $k=1$ general position means
that the slopes $r_{l,0}$ are pairwise distinct (\lean{injective\_of\_generalPosition\_one}); this is the
condition on the slopes in \cite{L}.

\begin{lemma}[\lean{exists\_hyperplane}]\label{lem:hyper}
Let $1\le s\le k$ and $u_1,\dots,u_s\in\Q^k$. There are $c\in\Q^k$, $c\ne0$, with $c_j=0$ for $j<k-s$,
and $c_*\in\Q$ such that $\sum_jc_ju_{x,j}+c_*=0$ for $x=1,\dots,s$.
\end{lemma}
\begin{proof}
The $s$ equations in the $s+1$ unknowns $c_{k-s},\dots,c_{k-1},c_*$ have a nonzero solution. If $c=0$
then $c_*=0$ by any of the equations.
\end{proof}

Fix for every $s$ and every $s$-tuple $u$ one such pair (\lean{hyper}), and let $\Phi_k(\mathcal L)$
(\lean{forbidden}) be the list of these pairs for all $1\le s\le k$ and all $s$-tuples of members of
$\mathcal L$. (The Lean list also contains the pair $(0,0)$, for $s=0$; it imposes no condition.)

\begin{lemma}[\lean{generalPosition\_succ}]\label{lem:gpsucc}
If $p_0,\dots,p_{n-1}\in\Q^k$ are in general position and $p_n$ is off $\Phi_k(p_{n-1},\dots,p_0)$, then
$p_0,\dots,p_n$ are in general position.
\end{lemma}
\begin{proof}
Let $S$, $t=|S|\le k+1$ and $\mu$ be as in Definition~\ref{def:gp}. If $n\notin S$ the hypothesis
applies. Let $n\in S$, $S_0=S\setminus\{n\}$ and $s=t-1\le k$. If $s=0$, then $\mu_n=0$ by the first
relation. If $s\ge1$, let $(c,c_*)$ be the pair of the $s$-tuple $(p_l)_{l\in S_0}$. Since $c_j=0$ for
$j<k-s=k+1-t$,
\[
\sum_{l\in S}\mu_l\Big(\sum_jc_jp_{l,j}+c_*\Big)=\sum_{j\ge k+1-t}c_j\sum_{l\in S}\mu_lp_{l,j}+c_*\sum_{l\in S}\mu_l=0 .
\]
The brackets with $l\in S_0$ vanish, and the bracket with $l=n$ does not; so $\mu_n=0$. The relations
then hold for $S_0$, with $t-1$ in place of $t$ (they are a subset of the given ones), and the
hypothesis gives $\mu_l=0$ on $S_0$.
\end{proof}

\begin{lemma}[Selection; \lean{exists\_pairs\_general\_position\_above\_thresholds}]\label{lem:sel}
Suppose a $\Phi_k$-supply of $\nu$-good pairs exists. Let $m\ge0$, let $X\in\R$, and let thresholds
$T_0,\dots,T_{m-1}$ be given, $T_i$ being a rational-valued function of $i$ rational arguments. Then
there are $\nu$-good pairs $(r_0,\varrho_0),\dots,(r_{m-1},\varrho_{m-1})$ whose slope vectors are in general
position, with
\[
q_i=q(r_i,\varrho_i)\ge2,\qquad w_i=\lceil\log q_i\rceil>X,\qquad w_i>T_i(w_0,\dots,w_{i-1})\qquad(0\le i<m).
\]
\end{lemma}
\begin{proof}
Choose the pairs one at a time (\lean{exists\_successive}, \lean{exists\_above\_thresholds}; these
statements do not depend on the nature of the objects selected): when the first $i$ pairs are chosen,
the supply gives a good pair whose slope vector is off $\Phi_k(r_{i-1},\dots,r_0)$ and whose denominator
is so large that the three inequalities hold. General position follows from Lemma~\ref{lem:gpsucc} by
induction (\lean{generalPosition\_of\_admissible}).
\end{proof}

The weights $w_i=\lceil\log q_i\rceil$ are positive integers with $w_i-1<\log q_i\le w_i$
(\lean{logWeights}).

\subsection{The statement that remains}

\begin{proposition}[\lean{CLogEndpointReduction.no\_supplyC}]\label{prop:no}
Let $k\ge1$ and assume $\Unif(k)$. Let $a_0,\dots,a_{k-1},g$ be positive rational numbers such that
$(1,(\log a_j),\log g)$ is independent, and $\nu>k+2$. Then there is no $\Phi_k$-supply of $\nu$-good
pairs.
\end{proposition}

Theorem~\ref{step} follows from Propositions~\ref{prop:nosupply} (with $\Phi=\Phi_k$) and~\ref{prop:no}.
Proposition~\ref{prop:no} is proved in Section~\ref{sec:par}, after the interpolation statement and the
two bounds for the determinant.

\section{The interpolation statement}\label{sec:int}

Files \lean{NLogIndex}, \lean{NLogInterpolation}, \lean{NLogMargins}, \lean{NLogCharacterCenters},
\lean{CLogInterpolation}, \lean{CLogTranslation}.

\subsection{Coordinates and the chart}
Let $k\ge0$ and $m\ge0$. The ambient affine space has the coordinates $Y_0,\dots,Y_k$ (torus
coordinates) and $X_0,\dots,X_{m-1}$ (additive coordinates). The index set is $\{0,\dots,m+k\}$; the
torus coordinate $j$ has index $j$ (\lean{torusIdx}), the additive coordinate $i$ has index $i+k+1$
(\lean{addIdx}). We write $x$ for a general index and say that it is a torus index or an additive
index. Let
\[
\mathcal P=\C[Y_0,\dots,Y_k,X_0,\dots,X_{m-1}],\qquad\mathcal S=\C[[t_0,\dots,t_k,u_0,\dots,u_{m-1}]]
\]
(\lean{Polynomial k m}, \lean{Series k m}). For slopes $r=(r_{i,j})$, $0\le i<m$, $0\le j<k$, put
\begin{equation}\label{forms}
\ell_i(z)=z_k-\sum_{j<k}r_{i,j}z_j\qquad(z\in\C^{k+1}).
\end{equation}
A \emph{center} is a point $c=(b,\xi)$ with the torus part $b\in\C^{k+1}$ and the additive part
$\xi\in\C^m$ (\lean{fullCenterC}); it is \emph{on the torus} if $\xi=0$ (\lean{fullCenter}). The
\emph{chart} at $c$ is the homomorphism $\mathcal P\to\mathcal S$ (\lean{formalJetC r b ξ}) with
\begin{equation}\label{chart}
Y_j\longmapsto b_j(1+t_j),\qquad
X_i\longmapsto\xi_i+u_i+\log(1+t_k)-\sum_{j<k}r_{i,j}\log(1+t_j)=\xi_i+u_i+\ell_i\big(\log(1+t)\big),
\end{equation}
where $\log(1+t)=(\log(1+t_0),\dots,\log(1+t_k))$. It is the chart at the center $(b,0)$ on the torus
(\lean{formalJet r b}) composed with the translation
\begin{equation}\label{tau}
\tau_\xi:\mathcal P\to\mathcal P,\qquad(\tau_\xi P)(Y,X)=P(Y,X+\xi)
\end{equation}
(\lean{shiftX}), which preserves the degree of a polynomial in every variable and commutes with all
partial derivatives.

\subsection{Numerical data}

\begin{definition}[\lean{Margins}]\label{def:margins}
For integers $k\ge1$, $m,K\ge0$ and rationals $\omega,v,\theta,\sigma$ the proposition
$\mathrm{Mar}_k(m,K)=\mathrm{Mar}_k(m,K;\omega,v,\theta,\sigma)$ is defined by recursion on $k$.
$\mathrm{Mar}_1(m,K)$ is the conjunction of
\begin{gather*}
2\le m,\qquad(1+3\sigma)^{m+2}K\Big(\frac{\omega}{v}\Big)^{2}\theta^{m}<1,\qquad
(1+3\sigma)^{m+1}\,\frac{\omega}{v}\,\theta^{m-1}<1,\\
(1+3\sigma)^{m-1}\theta^{m-1}<1,\qquad\frac{4(m+3)}{\sigma}<\frac{\omega}{v};
\end{gather*}
these are the margins of the two-torus lemma of \cite{L}. For $k\ge2$, $\mathrm{Mar}_k(m,K)$ is the
conjunction of
\begin{equation}\label{margins}
k+1\le m,\qquad(1+3\sigma)^{m+k+1}K\Big(\frac{\omega}{v}\Big)^{k+1}\theta^{m}<1,\qquad
\frac{(k+2)(m+k+2)}{\sigma}<\frac{\omega}{v},
\end{equation}
and $\mathrm{Mar}_{k-1}(m-1,1)$.
\end{definition}

The second condition in \eqref{margins} is the volume margin for $k+1$ torus coordinates and the third
is the torus margin. Unrolling the recursion, $\mathrm{Mar}_k(m,K)$ consists of the conditions
\eqref{margins} for the pairs $(k-s,m-s)$, $0\le s\le k-2$, with $K$ replaced by $1$ when $s\ge1$, and of
$\mathrm{Mar}_1(m-k+1,1)$.

\begin{definition}[\lean{BaseData k}, \lean{BaseData.fibreData}]\label{def:data}
\emph{Data} for $k$ slopes consist of integers $m\ge0$, $K\ge1$ and rationals $\omega,v,\sigma>0$,
$0<\theta<1$ with
\[
(1+\sigma)\theta<1\qquad\text{and}\qquad\mathrm{Mar}_k(m,K;\omega,v,\theta,\sigma).
\]
(For $k=0$ the margins are defined to be false, so data exist only for $k\ge1$.) In Lean the torus
weight $\omega$ is the field \lean{w0}; in this paper $w_0,\dots,w_{m-1}$ always denote the weights of the
additive coordinates. If
$D=(m,K,\omega,v,\theta,\sigma)$ are data for $k\ge2$ slopes, the \emph{fibre data} are
$D^{\rm fib}=(m-1,1,\omega,v,\theta,\sigma)$; they are data for $k-1$ slopes.
\end{definition}

\begin{lemma}[\lean{BaseData.torus\_margin}]\label{lem:torus}
For data $D$ for $k$ slopes,
\[
\frac{(k+1)(m+k+2)}{\sigma}<\frac{\omega}{v}.
\]
\end{lemma}
\begin{proof}
For $k=1$ the left side is $2(m+3)/\sigma<4(m+3)/\sigma$; for $k\ge2$ it is less than the left side of
the third condition in \eqref{margins}.
\end{proof}

In Lemma~\ref{lem:torus} the factor $k+1$ is the number of torus coordinates and $m+k+2$ is the number
of coordinates plus one. This is the form in which the torus margin is used (Step~5 of
Section~\ref{sec:curve}).

\subsection{Weights and packets}
Let $D$ be data and let $w_0,\dots,w_{m-1}$ be positive rational \emph{weights} of the additive
coordinates (\lean{w : Fin m → ℚ}). Put
\begin{equation}\label{weights}
W_x=\begin{cases}\omega&x\text{ a torus index}\\ w_i&x\text{ the additive index }i\end{cases}\qquad
V_x=\begin{cases}v&x\text{ a torus index}\\ w_i/\theta&x\text{ the additive index }i\end{cases}
\end{equation}
(\lean{sourceWeights}, \lean{NLogWeightSeparation.degreeWeights}; \lean{jetWeights}). The $W$-degree of
a monomial $Y^hX^\gamma$ is $\omega|h|+\sum_iw_i\gamma_i$, and the $V$-weight of $t^\lambda u^\beta$ is
$v|\lambda|+\sum_iw_i\beta_i/\theta$, where $|h|=\sum_jh_j$. Let $\mathcal P_{\le H}$ be the space of
polynomials all of whose monomials have $W$-degree at most $H$ (\lean{WeightedPolynomials}). The
\emph{packet} of a series is the family of its coefficients at the monomials of $V$-weight strictly
less than $H$ (\lean{rationalCoefficientPacket}). For centers $c_J=(b_J,\xi_J)$, $0\le J<K$, the
\emph{packet map} (\lean{packetMapC}; \lean{packetMap} for centers on the torus) sends
$P\in\mathcal P_{\le H}$ to the $K$ packets of its images under the charts \eqref{chart} at
$c_0,\dots,c_{K-1}$.

\subsection{Characters and thresholds}

\begin{definition}[\lean{characterValue}, \lean{CharacterDistinct}]\label{def:char}
For $r\in\Q^k$ let $N(r)=\prod_{t<k}\operatorname{den}(r_t)$ be the product of the reduced denominators
(\lean{charDen}) and $p_j(r)=\operatorname{num}(r_j)\cdot N(r)/\operatorname{den}(r_j)=r_jN(r)\in\Z$
(\lean{charExp}). The \emph{character} of $r$ is
\[
\chi_r=Y_k^{N(r)}\prod_{j<k}Y_j^{-p_j(r)} .
\]
Slopes $r_0,\dots,r_{m-1}\in\Q^k$ and centers $c_J=(b_J,\xi_J)$, $0\le J<K$, with
$b_J\in(\C^\times)^{k+1}$ have \emph{distinct characters} if for every $i$ the values $\chi_{r_i}(b_J)$,
$0\le J<K$, are pairwise distinct. This is a condition on the torus parts of the centers only.
\end{definition}

Formally $N(r)\,d\ell_r(\log Y)=d\log\chi_r$. The exponent $N(r)$ is the product, not the least common
multiple, of the denominators; this is immaterial for the argument.

\begin{lemma}[\lean{characterDistinct\_geometric}]\label{lem:centers}
Let $\beta_0,\dots,\beta_k$ be positive real numbers without a power relation and
$b_J=(\beta_0^J,\dots,\beta_k^J)$, $0\le J<K$. Then any rational slopes $r_0,\dots,r_{m-1}$ and any
centers with the torus parts $b_J$ have distinct characters.
\end{lemma}
\begin{proof}
$\chi_r(\beta)$ is a positive real number, and it is not $1$ because the exponent $N(r)$ of $\beta_k$ is
not zero (\lean{realCharacter\_ne\_one}). Since $\chi_r(b_J)=\chi_r(\beta)^J$, the values are distinct.
\end{proof}

\emph{Thresholds} for $m$ weights (\lean{SeparatedWeightPA.Thresholds}, from \cite{L}) are functions
$T_0,\dots,T_{m-1}$, where $T_i$ is a rational-valued function of $i$ rational arguments. Weights
$w_0,\dots,w_{m-1}$ are \emph{above} the thresholds (\lean{AboveThresholds}) if $w_i>0$ and
$w_i>T_i(w_0,\dots,w_{i-1})$ for all $i$.

\begin{definition}[\lean{UniformStatementC k}]\label{def:U}
$\Unif(k)$ is the following statement. For all data $D=(m,K,\omega,v,\theta,\sigma)$ for $k$ slopes there
are thresholds $T$ for $m$ weights such that for all rational weights $w$ above $T$ there is an integer
$R\ge1$ such that for all slopes $r_0,\dots,r_{m-1}\in\Q^k$ in general position and all centers
$c_J=(b_J,\xi_J)$, $0\le J<K$, with $b_J\in(\C^\times)^{k+1}$, arbitrary $\xi_J\in\C^m$ and distinct
characters there is $H_0$ such that the packet map
\[
\mathcal P_{\le H}\longrightarrow\bigoplus_{J<K}\ \big(\text{packets of $V$-weight}<H\big)
\]
is surjective for every integer $H\ge H_0$ divisible by $R$.
\end{definition}

The order of the quantifiers is: data; thresholds; weights; modulus $R$; slopes and centers; $H_0$.
The thresholds and the modulus do not depend on the slopes and the centers. For centers on the torus
the statement is \lean{UniformStatement k} (\lean{uniformStatement\_of\_C}); for $k=1$ and centers on
the torus it is the two-torus interpolation lemma of \cite{L}, which is stated there for nonzero
slopes.

\section{The curve inequality}\label{sec:curve}

\subsection{Separation of weights}
Files \lean{NLogWeightSeparation}, \lean{NLogDeletionSeparation}, \lean{NLogLowCodimensionExclusion},
\lean{CLogCurveBaseOne}, \lean{CLogCurveChain}, \lean{CLogGeometricData}.

\begin{definition}[\lean{SeparatedProducts D c w}]\label{def:sep}
Let $D$ be data for $k$ slopes, $c$ a rational number and $w$ positive weights. The weights are
\emph{$c$-separated} if for all sets $A,A'$ of indices with $|A|=|A'|$ and every additive index $x$
such that $x\in A$, $x\notin A'$, and every index larger than $x$ belongs to $A$ if and only if it
belongs to $A'$,
\[
c\prod_{y\in A'}V_y<\prod_{y\in A}W_y .
\]
\end{definition}

So the largest index at which $A$ and $A'$ differ is an additive index and lies in $A$. No condition is
imposed when $A$ and $A'$ differ in torus indices only: the $k+1$ torus weights are tied.

\begin{definition}[\lean{frameComparisonConstant}, \lean{baseConstantOne}, \lean{chainConstantOne},
\lean{chainConstantK}]\label{def:chain}
For data $D$ for $k$ slopes put
\[
F(D)=2\,P^{P}\Big(1+\frac P\sigma\Big)^{P},\qquad P=m+k+2 .
\]
The \emph{chain constant} $c_k(D)$ is defined for $k\ge1$ by
\[
c_1(D)=\max\big(F_2(D),\ F'(D)\big),\qquad
c_k(D)=\max\big(F(D),\ c_{k-1}(D^{\rm fib})\big)\quad(k\ge2),
\]
where $F_2(D)=2P^{P}(1+P/\sigma)^{P}$ with $P=m+3$ is the frame constant of the two-torus curve
inequality of \cite{L} (\lean{TwoLogNormalFrameRigidity.frameComparisonConstant}) and
$F'(D)=2P'^{P'}(1+P'/\sigma)^{P'}$ with $P'=m+1$ is the separation constant of the one-torus proposition
of \cite{L} for $m-1$ additive coordinates (\lean{PAWeightGeometry.separationConstant}). Finally
\[
c^+_k(D)=\max\big(F(D),\ c_k(D)\big)\qquad(k\ge1).
\]
\end{definition}

All the terms have the form $2P^P(1+P/\sigma)^P$, with $P$ the number of coordinates plus one at the
corresponding stage, so $c^+_k(D)=c_k(D)=F(D)$; the formal proof does not use this and works with the
maxima. The constant $c_k(D)$ is used for centers on the torus and $c^+_k(D)$ for arbitrary centers.

\begin{lemma}[\lean{CLogGeometry.GeometricData.exists\_separation\_thresholds}]\label{lem:thr}
Let $k\ge1$ and let $D$ be data for $k$ slopes. There are thresholds $T$ such that all weights above
$T$ are $c^+_k(D)$-separated.
\end{lemma}
\begin{proof}
As in \cite{L}: if the largest differing index is the additive index $i$, the ratio of the two sides
of the inequality of Definition~\ref{def:sep} is $w_i$ times a positive quantity which depends on the
data and on $w_0,\dots,w_{i-1}$ only; so the condition is one lower bound for each $w_i$ in terms of
the earlier weights (\lean{NLogWeightSeparation.exists\_separation\_thresholds}, for any positive
constant).
\end{proof}

\begin{lemma}[\lean{separatedProducts\_mono}, \lean{fibre\_separated}]\label{lem:del}
\textup{(a)} If the weights are $c$-separated and $c'\le c$, they are $c'$-separated.
\textup{(b)} Let $D$ be data for $k\ge2$ slopes, $w$ weights that are $c$-separated for $D$, and
$0\le i<m$. Then the weights $w'=(w_l)_{l\ne i}$ are $c$-separated for $D^{\rm fib}$.
\end{lemma}
\begin{proof}
(b) Send the torus index $j\le k-1$ of the fibre data to the torus index $j$ of $D$, and the additive
indices to the additive indices of $D$ other than $i$, in increasing order. This map of index sets is
strictly increasing (\lean{deletionIndex}) and preserves $W$ and $V$ (\lean{deletion\_weights}). The
inequality for $A,A'$ is the inequality for their images.
\end{proof}

\subsection{Degree and contact}
Let $E$ be a finitely generated field of transcendence degree $1$ over $\C$, and let
$u=(u_x)_{0\le x\le m+k}\in E^{m+k+1}$ be a tuple with at least one transcendental entry; write
$y_j=u_j$ for the torus entries and $x_i=u_{i+k+1}$ for the additive entries. Places $p$ of $E$ are the
points of its complete smooth model, with normalized valuations $\ord_p$ (\lean{NormalizedPlace}). Put
(\lean{CurveContactSum.weightedDegree})
\[
\deg_W(u)=\sum_p\max\Big(0,\ \max_x\Big(-\frac{\ord_pu_x}{W_x}\Big)\Big).
\]
A place $p$ is \emph{centered at} $c_J=(b_J,\xi_J)$ if all $u_x$ are regular at $p$ with
$y_j(p)=b_{J,j}$ and $x_i(p)=\xi_{J,i}$. Then $t_j=y_j/b_{J,j}-1$ lies in the maximal ideal of the completed local ring at $p$, and so
do the elements
\begin{equation}\label{branch}
u^{\circ}_i=x_i-\xi_{J,i}-\log(1+t_k)+\sum_{j<k}r_{i,j}\log(1+t_j)
\end{equation}
(\lean{branchCoordinates}); the \emph{contact} at $p$ (\lean{contactC}; \lean{NLogInterpolation.contact}
for centers on the torus) is
\[
h_p=\min\Big(\min_{j\le k}\frac{\ord_pt_j}{v},\ \min_{i<m}\frac{\ord_pu^\circ_i}{w_i/\theta}\Big),
\]
with $\ord_p0=+\infty$. At places not centered at any center put $h_p=0$. Only finitely many places are
centered (\lean{placesC}, \lean{places}).

\begin{proposition}[Curve inequality; \lean{weighted\_curve\_inequalityC\_chain},
\lean{weighted\_curve\_inequalityC\_of\_one\_le}]\label{prop:CI}
Let $k\ge1$, let $D$ be data for $k$ slopes and $w$ positive rational weights that are
$c^+_k(D)$-separated. Let $r_0,\dots,r_{m-1}\in\Q^k$ be in general position and let $c_J=(b_J,\xi_J)$,
$0\le J<K$, be centers with $b_J\in(\C^\times)^{k+1}$ and distinct characters. Then for every $E$ and $u$
as above
\[
(1+\sigma)\sum_ph_p\ \le\ \deg_W(u).
\]
\end{proposition}

\subsection{Reduced slopes}
File \lean{NLogReducedSlopes}.

\begin{definition}[\lean{fibreScale}, \lean{reducedSlope}]\label{def:red}
Let $k\ge1$, $r_0,\dots,r_{m-1}\in\Q^k$ and $0\le i<m$. For $l\ne i$ put
\[
d_l=r_{i,k-1}-r_{l,k-1},\qquad\rho_{l,j}=\frac{r_{l,j}-r_{i,j}}{d_l}\quad(0\le j<k-1),
\]
the second when $d_l\ne0$. The vectors $\rho_l\in\Q^{k-1}$, $l\ne i$, are the \emph{reduced slopes} with
respect to $i$.
\end{definition}

\begin{lemma}[\lean{fibreScale\_ne\_zero}, \lean{generalPosition\_reducedSlope}]\label{lem:heredity}
If $r_0,\dots,r_{m-1}\in\Q^k$ are in general position, then $d_l\ne0$ for all $l\ne i$, and the reduced
slopes $(\rho_l)_{l\ne i}$ are in general position in $\Q^{k-1}$.
\end{lemma}
\begin{proof}
General position for $S=\{i,l\}$ says that $(r_{i,k-1},1)$ and $(r_{l,k-1},1)$ are independent, so
$d_l\ne0$. Let $S'$ be a set of indices different from $i$ with $t'=|S'|\le k$, and $\nu_l$ $(l\in S')$
with $\sum\nu_l=0$ and $\sum\nu_l\rho_{l,j}=0$ for $k-t'\le j\le k-2$. Put
\[
\mu_l=\frac{\nu_l}{d_l}\ (l\in S'),\qquad\mu_i=-\sum_{l\in S'}\mu_l,\qquad S=S'\cup\{i\},\quad t=t'+1\le k+1 .
\]
Then $\sum_{l\in S}\mu_l=0$, and for every $j$
\[
\sum_{l\in S}\mu_lr_{l,j}=\sum_{l\in S'}\mu_l\,(r_{l,j}-r_{i,j}).
\]
For $j=k-1$ this is $-\sum_{l\in S'}\mu_ld_l=-\sum\nu_l=0$, and for $k-t'\le j\le k-2$ it is
$\sum_{l\in S'}\nu_l\rho_{l,j}=0$. These are the relations of Definition~\ref{def:gp} for $S$, since
$k+1-t=k-t'$. Hence $\mu=0$ on $S$ and $\nu=0$.
\end{proof}

\subsection{Excess contact and the character fibre}
The first seven steps of the proof are the same for every $k\ge1$ and for arbitrary centers.

\begin{proposition}[\lean{excess\_implies\_character\_fibreC}; \lean{excess\_implies\_character\_fibre}
for centers on the torus]\label{prop:excess}
Let $k\ge1$, let $D$ be data for $k$ slopes, $w$ positive rational weights that are $F(D)$-separated,
$r_0,\dots,r_{m-1}\in\Q^k$, and $c_J=(b_J,\xi_J)$, $0\le J<K$, centers with $b_J\in(\C^\times)^{k+1}$. Let $E$
and $u$ be as above, and suppose that $\deg_W(u)<(1+\sigma)\sum_ph_p$. Then there are a center $c_J$ and
an additive index $i$ such that $y_j\ne0$ for all $j$ and
\begin{equation}\label{fibre}
x_i=\xi_{J,i},\qquad\chi_{r_i}(y)=\chi_{r_i}(b_J)\in\C .
\end{equation}
\end{proposition}

Neither general position nor the distinctness of the characters is used here.

\begin{proof}
By the assumption some place is centered at a center.

\emph{Step 1: the frame} (\lean{NLogDifferentialFrame}, \lean{NLogFrameDegree}). The derivations of
$\mathcal P$
\[
D_j=Y_j\partial_{Y_j}-\sum_ir_{i,j}\partial_{X_i}\ (j<k),\qquad
D_k=Y_k\partial_{Y_k}+\sum_i\partial_{X_i},\qquad D_{i+k+1}=\partial_{X_i}\ (i<m)
\]
commute and do not increase the $W$-degree. Under the chart \eqref{chart} at any center they
correspond to $(1+t_j)\partial_{t_j}$ and $\partial_{u_i}$; they commute with the translations
\eqref{tau} (\lean{polynomialFrameWord\_shiftX}). Give $D_x$ the cost $V_x$; the cost of a
word in the $D_x$ is the sum of the costs of its letters. A word of cost $c$ lowers the $V$-order of
the image of a polynomial under a chart by at most $c$.

\emph{Step 2: auxiliary polynomial} (\lean{NLogAuxiliaryPolynomial}, \lean{CLogAuxiliary}). The number of monomials of
$W$-degree at most $N$ is asymptotic to $N^{m+k+1}/((m+k+1)!\,\omega^{k+1}\prod_iw_i)$, and the number of
monomials $t^\lambda u^\beta$ of $V$-weight less than $(1+3\sigma)N$, at $K$ centers, to
$K(1+3\sigma)^{m+k+1}N^{m+k+1}\theta^m/((m+k+1)!\,v^{k+1}\prod_iw_i)$. By the volume margin in
\eqref{margins}, for all large $N$ there is a nonzero $F\in\mathcal P_{\le N}$ whose image under the
chart at every center has $V$-order at least $(1+3\sigma)N$. (The conditions on $F$ at the center
$(b_J,\xi_J)$ are the conditions on $\tau_{\xi_J}F$ at $(b_J,0)$, and $\tau_{\xi_J}$ is an automorphism of
$\mathcal P_{\le N}$; the count does not depend on the additive parts of the centers.)

\emph{Step 3: vanishing} (\lean{NLogDerivativeVanishing}, \lean{CLogContact}). Let $G$ be obtained from $F$ by a word of
cost at most $\sigma N$. Then $G\in\mathcal P_{\le N}$ and its images have $V$-order at least
$(1+2\sigma)N$. At a place centered at $c_J$ the substitution of the branch \eqref{branch} into the
chart gives $\ord_pG(u)\ge(1+2\sigma)Nh_p$ (\lean{word\_order\_lowerC}); at any place
$\ord_pG(u)\ge-N\max(0,\max_x(-\ord_pu_x/W_x))$. If $G(u)\ne0$, the sum of its orders over all places is
zero, and $(1+2\sigma)N\sum h_p\le N\deg_W(u)<(1+\sigma)N\sum h_p$, which is impossible. So $G(u)=0$.

\emph{Step 4: a persistent component and the comparison} (\lean{NLogPersistentComponents},
\lean{NLogTaylorRectangular}, \lean{NLogFrameTransverseMultiplicity},
\lean{NLogPrimeFrameMultiplicity}, \lean{NLogTransverseSliceComparison},
\lean{NLogPersistentComparison}). Let $\mathfrak p\subset\mathcal P$ be the kernel of $P\mapsto P(u)$;
it is a prime ideal of height $m+k$. Put
\[
\epsilon=\frac{\sigma}{m+k+2},
\]
and let $I_s$, $0\le s\le m+k+2$, be the ideal generated by the images of $F$ under all words of cost at
most $s\,\epsilon N$. By Step~3, $I_s\subseteq\mathfrak p$ for all $s$. Since $F\ne0$ and the heights of
primes contained in $\mathfrak p$ lie between $1$ and $m+k$, there are $s<m+k+2$ and a prime
$\mathfrak q\subseteq\mathfrak p$ which is a minimal prime of both $I_s$ and $I_{s+1}$
(\lean{exists\_persistent\_frame\_component}). Every word of cost at most $\epsilon N$ maps $I_s$ into
$\mathfrak q$. Let $d$ be the height of $\mathfrak q$, $1\le d\le m+k$. A set $A$ of $d$ indices is a
\emph{coordinate normal basis} if the coordinate vectors $e_x$, $x\in A$, form a basis of the normal
space of $\mathfrak q$ (the quotient of the ambient space over the residue field of $\mathfrak q$ by the
tangent space of $\mathfrak q$), and a set $B$ is a \emph{frame normal basis} if the frame vectors
$D_x$, $x\in B$, do. The weighted transverse multiplicity lemma of \cite{L}, which holds for an
arbitrary commuting frame, gives for all large $N$ and all such $A$, $B$
\begin{equation}\label{mult}
\prod_{x\in A}W_x\ \le\ \frac{d^{\,d}}{\epsilon^{d}}\prod_{x\in B}V_x
\end{equation}
(\lean{NormalProductComparison}, \lean{exists\_persistent\_component\_with\_comparison}). Moreover
\begin{equation}\label{raw}
\frac{d^{\,d}}{\epsilon^{d}}=\Big(\frac{d\,(m+k+2)}{\sigma}\Big)^{d}\le F(D)\qquad(d\le m+k+1)
\end{equation}
(\lean{raw\_factor\_le\_frameComparisonConstant}).

\emph{Step 5: the additive directions are not all tangent} (\lean{NLogLowCodimensionExclusion}).
Suppose every $\partial_{X_i}$ is tangent to $\mathfrak q$. Then the normal space is spanned by the
images of $\partial_{Y_0},\dots,\partial_{Y_k}$, so $d\le k+1$, and modulo the tangent space $D_j$ is the
nonzero multiple $Y_j\partial_{Y_j}$ of $\partial_{Y_j}$ for every torus index $j$. A set of $d$ torus
indices whose directions form a basis of the normal space is both a coordinate and a frame normal basis,
and \eqref{mult} reads $(\omega/v)^d\le(d(m+k+2)/\sigma)^d$. Hence
$\omega/v\le(k+1)(m+k+2)/\sigma$, against Lemma~\ref{lem:torus}. With the factor $4$ of \cite{L}
in place of $k+1$ this argument works only for $d\le4$, that is, for at most four torus coordinates;
this is the reason for the factor in Definition~\ref{def:margins}.

\emph{Step 6: a logarithmic differential} (\lean{NLogNormalFrameRigidity}). Let $i$ be the largest
additive coordinate with $\partial_{X_i}$ not tangent. Its image extends to a coordinate normal basis
$A$ containing the index of $i$; larger additive indices occur in no normal basis of either kind. If a
frame normal basis $B$ omitted the index of $i$, the weights being $F(D)$-separated would contradict
\eqref{mult} and \eqref{raw}. So every frame normal basis contains it, the other frame vectors do not
span the normal space, and the tangent space of $\mathfrak q$ lies in their ambient span, which is the
kernel of
\[
dX_i-\frac{dY_k}{Y_k}+\sum_{j<k}r_{i,j}\frac{dY_j}{Y_j}
\]
(\lean{logarithmicForm}). With the integers of Definition~\ref{def:char},
$N(r_i)\,dX_i=d\log\chi_{r_i}$ on $\mathfrak q$.

\emph{Step 7: rigidity} (\lean{NLogCharacterRigidity}, \lean{NLogPersistentRigidity},
\lean{NLogCurveRigidity}). As in \cite{L}: an exact differential has no residues, a logarithmic
differential of a nonconstant function has one; so $\chi_{r_i}$ and $X_i$ are constant modulo
$\mathfrak q$, hence on $u$. Evaluating at a place centered at some $c_J$: all $y_j\ne0$, and
\eqref{fibre} holds.
\end{proof}

\begin{corollary}[\lean{character\_value\_at\_center}]\label{cor:one}
If, in Proposition~\ref{prop:excess}, the slopes and the centers have distinct characters, then every
centered place of $u$ is centered at the center $c_J$ of \eqref{fibre}.
\end{corollary}
\begin{proof}
At a place centered at $c_{J'}$ the constant $\chi_{r_i}(y)$ takes the value $\chi_{r_i}(b_{J'})$; so
$\chi_{r_i}(b_{J'})=\chi_{r_i}(b_J)$ and $J'=J$.
\end{proof}

\subsection{Centers on the torus: induction on the number of torus coordinates}
In this subsection all centers are on the torus, $c_J=(b_J,0)$, and the weights are
$c_k(D)$-separated. Proposition~\ref{prop:CI} for such centers and weights is proved by induction on
$k\ge1$ (\lean{weighted\_curve\_inequality\_chain\_aux}).

\emph{The case $k=1$} is the curve inequality for two torus coordinates of \cite{L}
(\lean{TwoLogInterpolation.weighted\_curve\_inequality}). For $k=1$ the definitions of this section
are those of \cite{L} (file \lean{CLogCurveBaseOne}): the data are the data of the two-torus lemma,
general position is the injectivity of the slopes, the character of a slope $r$ is
$Y_1^{\operatorname{den}r}Y_0^{-\operatorname{num}r}$, the charts, the elements \eqref{branch} and the
contacts coincide, and $c_1(D)$-separated weights are separated with each of the two constants used
there (\lean{curve\_inequality\_one}).

\emph{The inductive step.} Let $k\ge2$, assume the statement for $k-1$, and suppose
$\deg_W(u)<(1+\sigma)\sum_ph_p$. By Lemma~\ref{lem:del}(a) the weights are $F(D)$-separated and
$c_{k-1}(D^{\rm fib})$-separated. By Proposition~\ref{prop:excess} and Corollary~\ref{cor:one} there are
$J$ and $i$ with \eqref{fibre}, where now $x_i=0$, and every centered place of $u$ is centered at
$b_J$.

\emph{Step 8: reduction to $k$ torus coordinates} (\lean{NLogCharacterFiberField},
\lean{NLogCharacterFiberBranch}, \lean{NLogCharacterFiberPlace}, \lean{NLogCharacterFiberPole},
\lean{NLogCharacterFiberSum}). Let $d_l$ and $\rho_l$ $(l\ne i)$ be as in Definition~\ref{def:red};
$d_l\ne0$ by Lemma~\ref{lem:heredity}. Consider the \emph{projected tuple}
(\lean{projectedCoordinates})
\[
u'=\Big(y_0,\dots,y_{k-1},\ \big(x_l/d_l\big)_{l\ne i}\Big)\in E^{(m-1)+(k-1)+1},
\]
obtained by deleting the last torus coordinate and the additive coordinate $i$ and rescaling. We
apply the induction hypothesis to $u'$ with the data $D^{\rm fib}$, the weights $w'=(w_l)_{l\ne i}$, the
slopes $\rho_l$ and the single center $b'=(b_{J,0},\dots,b_{J,k-1})$.
\begin{itemize}[leftmargin=2em]
\item \emph{Hypotheses.} $w'$ is $c_{k-1}(D^{\rm fib})$-separated for $D^{\rm fib}$ by
  Lemma~\ref{lem:del}(b); the $\rho_l$ are in general position by Lemma~\ref{lem:heredity}; with one
  center the condition on characters is empty.
\item \emph{$u'$ has a transcendental entry} (\lean{projected\_transcendental}). Otherwise
  $y_0,\dots,y_{k-1}$ and all $x_l$ $(l\ne i)$ are constants; $x_i=0$; and $y_k$ is algebraic over $\C$ by
  \eqref{fibre}, since $N(r_i)\ne0$. Then $u$ has no transcendental entry.
\item \emph{Contacts} (\lean{contact\_sum\_le\_projected\_at\_center}). At a place $p$ centered at $b_J$,
  \eqref{fibre} reads $(1+t_k)^{N}=\prod_{j<k}(1+t_j)^{p_j}$ in the completed local ring, with $N=N(r_i)$,
  $p_j=p_j(r_i)$. Taking logarithms of $1$-units (\lean{logOf\_integer\_cleared}),
  \[
  \log(1+t_k)=\sum_{j<k}r_{i,j}\log(1+t_j).
  \]
  Substituting in \eqref{branch}, for $l\ne i$ (\lean{affine\_branch\_eq\_scaled})
  \[
  u^\circ_l=x_l-\sum_{j<k}(r_{i,j}-r_{l,j})\log(1+t_j)
  =d_l\Big(\frac{x_l}{d_l}-\log(1+t_{k-1})+\sum_{j<k-1}\rho_{l,j}\log(1+t_j)\Big).
  \]
  The bracket is the element \eqref{branch} for the tuple $u'$, the slopes $\rho$ and the center $b'$. So
  the contact $h'_p$ of $u'$ is the minimum over a subfamily of the quantities whose minimum is $h_p$:
  $h_p\le h'_p$. Every centered place of $u$ is centered at $b_J$, hence centered at $b'$ for $u'$, and
  $h'\ge0$ everywhere; so $\sum_ph_p\le\sum_ph'_p$.
\item \emph{Degrees} (\lean{weightedDegree\_projected\_le}). $\deg_{W'}(u')\le\deg_W(u)$: the maximum is
  over fewer coordinates and $\ord_p(x_l/d_l)=\ord_px_l$.
\end{itemize}
By the induction hypothesis,
\[
(1+\sigma)\sum_ph_p\le(1+\sigma)\sum_ph'_p\le\deg_{W'}(u')\le\deg_W(u),
\]
which contradicts the assumption. \qed

\begin{remark}
The function field $E$ and its places are the same at all stages of the induction; only the tuple of
coordinates changes. After $k-1$ reductions the tuple has two torus coordinates and the inequality
of \cite{L} applies, whose proof makes one more reduction of the same kind and ends with the
one-torus proposition of \cite{L}. The successive data are $D,D^{\rm fib},(D^{\rm fib})^{\rm fib},\dots$;
this is why the margins of Definition~\ref{def:margins} and the constant of
Definition~\ref{def:chain} are defined by recursion.
\end{remark}

\subsection{Centers with additive coordinates}
Files \lean{CLogContact}, \lean{CLogShiftDegree}, \lean{CLogCurveInequality}, \lean{CLogCurveChain}.
Now let the centers $c_J=(b_J,\xi_J)$ be arbitrary, as in Proposition~\ref{prop:CI}, and suppose
$\deg_W(u)<(1+\sigma)\sum_ph_p$ (\lean{weighted\_curve\_inequalityC\_chain\_aux}). The weights are
$F(D)$-separated and $c_k(D)$-separated. By Proposition~\ref{prop:excess} and
Corollary~\ref{cor:one}, every centered place of $u$ is centered at one center $c_J$.

\emph{Step 9: translation.} Consider the tuple (\lean{shiftTuple})
\[
u^{T}=\big(y_0,\dots,y_k,\ x_0-\xi_{J,0},\dots,x_{m-1}-\xi_{J,m-1}\big).
\]
It has a transcendental entry (\lean{shiftTuple\_transcendental}). A place of $u$ centered at $c_J$ is
a place of $u^T$ centered at $(b_J,0)$, with the same elements \eqref{branch} and the same contact.
Let $h^T_p$ be the contacts of $u^T$ with respect to the centers $(b_0,0),\dots,(b_{K-1},0)$ on the
torus; they are nonnegative, so $\sum_ph_p\le\sum_ph^T_p$ (\lean{contactC\_sum\_le\_shift}). Subtracting
constants does not increase the poles, $\ord_p(x-c)\ge\min(\ord_px,0)$, so
$\deg_W(u^T)\le\deg_W(u)$ (\lean{weightedDegree\_shiftTuple\_le}). The inequality for centers on the
torus, applied to $u^T$ with the same data, weights and slopes and the centers $(b_J,0)$, $0\le J<K$,
gives
\[
(1+\sigma)\sum_ph_p\le(1+\sigma)\sum_ph^T_p\le\deg_{W}(u^T)\le\deg_W(u),
\]
which contradicts the assumption. This proves Proposition~\ref{prop:CI}. \qed

\section{From the curve inequality to surjectivity}\label{sec:packets}

Theorem~\ref{thm:U} follows from Proposition~\ref{prop:CI} by the argument of
\cite[Section ``From curves to all coefficient packets'']{L}: embed the affine space by
all monomials of $W$-degree at most $R$, blow up an ideal sheaf supported at the centers, deduce from
the curve inequality that $A-(1+\sigma)E$ is nef, show that $aA-E$ is ample for some $a$, conclude that
$A-E$ is ample, and finish with Serre vanishing and the Rees algebra. The thresholds of $\Unif(k)$ are
those of Lemma~\ref{lem:thr} (\lean{exists\_numeric\_thresholds}). The modulus $R$ is a positive integer
such that $R/W_x$ and $R/V_x$ are integers for all indices $x$ (\lean{NumericData.scale}); it depends on
the data and the weights only. The steps that refer to the coordinates are the following; they are
formalized for all $k\ge1$ and for centers with additive coordinates. The file names quoted below
are those of the layer for centers on the torus; the files for arbitrary centers have the same
names with the prefix \lean{CLog}, where they exist (Appendix~\ref{app:lean}).
\begin{enumerate}[label=(\roman*),leftmargin=2.4em]
\item \emph{Algebraic local coordinates} (\lean{NLogPolynomialChart}, \lean{NLogJetIdeals}). For each
  $i$ choose one integer $T_i$ with $w_i/\theta<T_i\,v$ (\lean{logCutoff}); the same cutoff is used for
  all $k+1$ torus coordinates. Let $G_i$ be the Taylor polynomial of $\log(1+t)$ of degree less than
  $T_i$ and
  \[
  u'_i=X_i-\xi_{J,i}-G_i(t_k)+\sum_{j<k}r_{i,j}G_i(t_j),\qquad t_j=Y_j/b_{J,j}-1 .
  \]
  Then $u'_i-u_i$ is a sum of $k+1$ series in one variable each, all of whose monomials have $V$-weight
  greater than $w_i/\theta$. The substitution $u_i\mapsto u'_i$ and its inverse preserve every
  filtration ideal $F_V^{\ge H}$ (spanned by the monomials of $V$-weight at least $H$).
\item \emph{Contacts} (\lean{NLogTruncatedBranchContact}, \lean{NLogBranchColength},
  \lean{NLogJetColength}). On a branch, each tail has order at least $T_i\min_j\ord t_j$, so
  $\ord(u'_i-u^\circ_i)/(w_i/\theta)>\min_j\ord t_j/v$ or the tail vanishes; hence $h_p$ may be computed
  with the $u'_i$.
\item \emph{The ideal} (\lean{NLogCompactJetIdeal}, \lean{NLogJetLocalization},
  \lean{NLogAlgebraicJetPackets}). $I_J=\big(t_j^{R/v}\ (j\le k),\ (u'_i)^{R\theta/w_i}\ (i<m)\big)$ has
  finite colength; its valuation on a branch through $c_J$ is $Rh_p$; every generator has $V$-weight
  $R$, so $I_J^{\,n}\subseteq F_V^{\ge nR}$. The centers are distinct points, because their torus parts are: centers with
  distinct characters have distinct torus parts when $m\ge1$ (\lean{centers\_injective},
  \lean{fullCenterC\_injective}). So the $I_J$ define a coherent ideal sheaf.
\item \emph{Curves} (\lean{NLogExceptionalDegree}, \lean{NLogCurveImage}, \lean{NLogCurveModel},
  \lean{NLogBlowupAmple}). A curve in the projective closure that misses the affine chart misses the
  exceptional divisor. A curve meeting the chart gives a tuple $u$ in its function field with a
  transcendental entry, to which Proposition~\ref{prop:CI} applies.
\end{enumerate}
The coordinates $t_j,u'_i$, the ideals $I_J$ and the local rings at a center $(b_J,\xi_J)$ are the
images under the translation $\tau_{-\xi_J}$ of those at the center $(b_J,0)$ (\lean{coordinatesC},
\lean{powerIdealC}, \lean{map\_powerIdealC\_shiftX}). In this way the statements of (i)--(iii) and the
formula for the degree of the exceptional divisor on a curve in (iv) are deduced from the statements
for centers on the torus (files \lean{CLogJetIdeals}, \lean{CLogJetLocalization},
\lean{CLogJetColength}, \lean{CLogCompactJetIdeal}, \lean{CLogCenterIdentification},
\lean{CLogExceptionalDegree}). The global part of the argument (the blow-up, the degrees on curves,
ampleness, sections and surjectivity) is formalized for a record of geometric data which contains
the additive parts of the centers (\lean{CLogGeometricData} and 22 further files,
Appendix~\ref{app:lean}). The result is \lean{CLogMatrixTheorems} and
\lean{CLogUniformInterpolation.cLogUniformInterpolation}. \qed

\section{The determinant and its arithmetic lower bound}\label{sec:ar}

Files \lean{CLogMatrix}, \lean{CLogMatrixJet}, \lean{CLogJetDeterminant}, \lean{CLogArithmetic},
\lean{CLogArithmeticBounds}, \lean{CLogData}, \lean{CLogDataCounting}, \lean{CLogActualMinor},
\lean{CLogDataMinors}; \lean{NLogMatrix}, \lean{NLogMatrixShear}, \lean{NLogShearFiltration},
\lean{NLogMatrixJet}, \lean{NLogJetDeterminant}, \lean{NLogCounting}, \lean{NLogArithmetic},
\lean{NLogArithmeticBounds}.

\subsection{The setting of Sections \ref{sec:ar}--\ref{sec:par}}\label{sec:setting}
Let $k\ge1$. The following objects are fixed (structure \lean{CLogData k nu Lambda c}).
\begin{itemize}[leftmargin=2em]
\item Positive rational numbers $\beta_0,\dots,\beta_k$ (\lean{bases}) with logarithms
  $L_t=\log\beta_t$ (\lean{logs}), and $L_0\ne0$. In the application $\beta_j=a_j$ for $j<k$ and
  $\beta_k=g$.
\item A real number $\nu$, parameters $\theta,A,B,C,\eta$ as in Definition~\ref{def:par} (\lean{base}),
  data $D=(m,K,\omega,v,\theta,\sigma)$ for $k$ slopes with the same $\theta$ (\lean{geometry}), and a real
  number $\eps_0$ with $0<\eps_0<g$, $\eps_0\le\frac12$, where $g$ is the gap \eqref{gap}.
\item Pairs $(r_i,\varrho_i)\in\Q^k\times\Q$, $0\le i<m$, which are $\nu$-good for
  $\big((L_j)_{j<k},L_k\big)$, with the slope vectors $r_0,\dots,r_{m-1}$ in general position and
  $q_i=q(r_i,\varrho_i)\ge2$ (\lean{r}, \lean{rho}, \lean{q}). The weights are
  $w_i=\lceil\log q_i\rceil$ (\lean{rationalWeights}), and $w_*>0$ is a number with $w_*\le w_i$ for all
  $i$ (\lean{wstar}).
\end{itemize}
The forms are $\ell_i$ of \eqref{forms}; so $|r_{i,j}|\le1$ and $|\ell_i(L)-\varrho_i|\le q_i^{-\nu}$. The
centers are
\begin{equation}\label{centers}
c_J=\big(\beta_0^J,\dots,\beta_k^J;\ J\varrho_0,\dots,J\varrho_{m-1}\big),\qquad0\le J<K
\end{equation}
(\lean{xi}): the multiples of the rational point $(\beta;\varrho)$ of the group $\Gm^{k+1}\times\Ga^m$. We
write $b_J=(\beta_0^J,\dots,\beta_k^J)$ for their torus parts. Three numerical margins, which are also fields of
the structure, are stated in Section~\ref{sec:par}.

\subsection{The matrix}
Put
\[
F_0=\frac4\theta,\qquad T_i=\Big\lceil\frac{F_0w_i}{v}\Big\rceil,\qquad
G_i(x)=\sum_{1\le n<T_i}(-1)^{n+1}\frac{x^n}{n},\qquad\mathrm{L}_i=\operatorname{lcm}(1,\dots,T_i-1)
\]
(\lean{F0}, \lean{truncationOrders}, \lean{truncatedLog}, \lean{Nat.lcmUpto}). For real $H>0$ the
\emph{columns} are the monomials $Y^hX^\gamma$, $h\in\Z_{\ge0}^{k+1}$, $\gamma\in\Z_{\ge0}^m$, with
$\omega|h|+\sum_iw_i\gamma_i\le H$ (\lean{Column}); the \emph{rows} are the triples $(J,\lambda,\beta)$ with
$0\le J<K$, $\lambda\in\Z_{\ge0}^{k+1}$, $\beta\in\Z_{\ge0}^m$ and $v|\lambda|+\sum_iw_i\beta_i/\theta<H$
(\lean{Row}). The matrix $A_H^{\rm tr}$ (\lean{CLogMatrix.matrixC}, \lean{entryC},
\lean{matrixC\_eq\_entryC}) has the entries
\begin{equation}\label{matrix}
[t^\lambda u^\beta]\ \prod_{j\le k}\big(\beta_j^J(1+t_j)\big)^{h_j}\ \prod_{i<m}
\Big(J\varrho_i+u_i+G_i(t_k)-\sum_{j<k}r_{i,j}G_i(t_j)\Big)^{\gamma_i},
\end{equation}
and $A_H^{\rm full}$ is the matrix with $\log(1+t_j)$ in place of $G_i(t_j)$, that is, the matrix of the
packet map of Section~\ref{sec:int} for the centers \eqref{centers} in the monomial bases
(\lean{fullLogMatrixC}). Let $M=M_H$ be the number of rows. A \emph{selection} is an
injective map from rows to columns; it defines an $M\times M$ minor of either matrix
(\lean{selectedMinorC}).

\begin{lemma}[The same minors; \lean{det\_fullLogC\_selected\_eq\_selectedMinorC}]\label{lem:shear}
For every selection the corresponding minors of $A_H^{\rm tr}$ and $A_H^{\rm full}$ are equal.
\end{lemma}
\begin{proof}
We have $w_i/\theta\le T_iv$ (\lean{data\_truncation\_weight}). Put $R_i(x)=G_i(x)-\log(1+x)$, a series
of order at least $T_i$. The truncated substitution is the full one followed by
$u_i\mapsto u_i+R_i(t_k)-\sum_{j<k}r_{i,j}R_i(t_j)$, which is the composition of $k+1$ substitutions,
each adding to every $u_i$ a series in one torus variable (\lean{NLogMatrixShear},
\lean{shiftSingle\_mem}). Every monomial of an added series has $V$-weight at least $vT_i\ge w_i/\theta$,
so each of these substitutions preserves the filtration by $V$-weight and acts on the packet of
weights less than $H$ by a unipotent triangular matrix, independent of the column
(\lean{NLogShearFiltration}, \lean{det\_shift\_restrict}). This operation acts on the packet of the
image of an arbitrary polynomial under the chart at a center on the torus; in the row of the center
$c_J$ it is applied to $\tau_{J\varrho}P$ (\lean{det\_selected\_fullLogC\_eq\_truncatedLog}).
\end{proof}

\begin{lemma}[\lean{CLogMinors.cofinal\_data\_minors}]\label{lem:minors}
Assume $\Unif(k)$, and for all data $D$ let $T(D)$ be thresholds as in Definition~\ref{def:U}. In the
setting of \ref{sec:setting}, if the weights $w_i$ are above $T(D)$, then for every $H_1$ there are an
integer $H\ge H_1$, $H>0$, and a selection whose minor $\Delta_H$ of $A_H^{\rm tr}$ is not zero.
\end{lemma}
\begin{proof}
The $\beta_t$ are positive with independent logarithms, so the torus parts $b_J$ of the centers are
nonzero and the centers have distinct characters with any rational slopes (Lemma~\ref{lem:centers});
the slopes are in general position. By $\Unif(k)$, applied to the centers \eqref{centers}, the packet
map is surjective for all large $H$ divisible by $R$, so $A_H^{\rm full}$ has rank $M$ and a nonzero
minor (\lean{exists\_minorC\_of\_packets}); apply
Lemma~\ref{lem:shear}.
\end{proof}

The function $D\mapsto T(D)$ is fixed before the approximation pairs are chosen.

\subsection{Counting}
By counting lattice points in weighted simplices (\lean{tendsto\_rowCount\_normalized},
\lean{meanRowWeight\_le\_theta}),
\begin{equation}\label{count}
M_H\sim\frac{K\theta^mH^{m+k+1}}{(m+k+1)!\,v^{k+1}\prod_iw_i}\quad(H\to\infty),\qquad
\overline b_H:=\frac1{M_HH}\sum_{\rm rows}\sum_iw_i\beta_i\in[0,\theta]
\end{equation}
(\lean{actualRowCount}, \lean{actualMean}).

\subsection{The arithmetic lower bound}
Let $\Lambda_0=4+\log4$ (\lean{Arithmetic.lcmConstant}); then
$\log\operatorname{lcm}(1,\dots,n)\le\Lambda_0n$ for all $n$. Put
\begin{equation}\label{Ear}
E_{\rm ar}=\frac{\Lambda_0F_0}{v}+\frac{(K-1)\sum_{t\le k}\log\operatorname{den}(\beta_t)}{\omega}
\end{equation}
(\lean{arithmeticError}).

\begin{proposition}[\lean{CLogComparison.actual\_minor\_lower\_bound}]\label{prop:lower}
For every $H>0$ and every selection with $\Delta_H\ne0$,
\[
\frac{\log|\Delta_H|}{M_HH}\ \ge\ -(1-\overline b_H)-E_{\rm ar}.
\]
\end{proposition}
\begin{proof}
An entry \eqref{matrix} is zero unless $\beta_i\le\gamma_i$ for all $i$
(\lean{entryC\_eq\_zero\_of\_not\_le}); a nonzero term of the determinant expansion pairs every row with a
column with this property, so
\[
J_i:=\sum_{\rm selected\ columns}\gamma_i-\sum_{\rm rows}\beta_i\ \ge\ 0
\]
(\lean{selectedMinorC\_total\_row\_degree\_le\_total\_column\_degree}). Since $q_ir_{i,j}\in\Z$ and
$q_i\varrho_i\in\Z$, the polynomial
$q_i\mathrm{L}_i\big(J\varrho_i+G_i(t_k)-\sum_{j<k}r_{i,j}G_i(t_j)\big)$ has integer coefficients
(\lean{cleared\_formC\_mem}, \lean{q\_mul\_xi\_mem}). Let $\mathrm D_H=\big(\prod_t\operatorname{den}(\beta_t)\big)^{(K-1)\lfloor H/\omega\rfloor}$
(\lean{baseDenominator}); since $J\le K-1$ and $h_j\le H/\omega$, the number
$\mathrm D_H\prod_j\beta_j^{Jh_j}$ is an integer. Multiplying the column of $X^\gamma$ by
$\prod_i(q_i\mathrm L_i)^{\gamma_i}$, the row of $\beta$ by $\prod_i(q_i\mathrm L_i)^{-\beta_i}$ and every
entry by $\mathrm D_H$ makes all entries integers (\lean{selectedMinorC\_entries\_gaussian}; in Lean the
cleared entries are shown to be Gaussian integers, which is what the general lemma of \cite{L}
provides and is sufficient). Hence
\[
\mathrm D_H^{\,M}\prod_i(q_i\mathrm L_i)^{J_i}\cdot|\Delta_H|\ \ge\ 1
\]
(\lean{selectedMinorC\_clearing\_bound}). Now $\log q_i\le w_i$,
$\log\mathrm L_i\le\Lambda_0(T_i-1)\le\Lambda_0F_0w_i/v$, and $\sum_iw_iJ_i\le MH-\sum_{\rm rows}\sum_iw_i\beta_i
=MH(1-\overline b_H)$ because every column has $\sum_iw_i\gamma_i\le H$ (\lean{combined\_log\_costC\_le});
and $\log\mathrm D_H\le(K-1)\sum_t\log\operatorname{den}(\beta_t)\,H/\omega$
(\lean{log\_baseDenominator\_le}).
\end{proof}

The constants $\varrho_i$ enter the arithmetic only through $q_i\varrho_i\in\Z$: this is why the
denominator of a pair is the common denominator of its slopes and of its constant.

\section{The analytic upper bound}\label{sec:an}

\subsection{Exact forms}
Files \lean{CLogAnalyticDataIndices}, \lean{CLogSlopeMatrix}. Every approximation pair gets its own
linear form which takes the value $\varrho_i$ at the point $L=(L_0,\dots,L_k)$. Put (\lean{formValue}, \lean{eta},
\lean{exactCoeffR})
\begin{equation}\label{eta}
\eta_i=\frac{\ell_i(L)-\varrho_i}{L_0},\qquad\ell_i^*(z)=\ell_i(z)-\eta_iz_0=\sum_{t\le k}c^*_{i,t}z_t ,
\end{equation}
so that $c^*_{i,k}=1$, $c^*_{i,j}=-r_{i,j}$ for $0<j<k$ and $c^*_{i,0}=-r_{i,0}-\eta_i$ (if $k=1$ the
index $0$ is the only slope index). Then
\begin{equation}\label{exact}
\ell_i^*(L)=\varrho_i,\qquad|\eta_i|\le\frac{q_i^{-\nu}}{|L_0|}\le\frac1{|L_0|},\qquad
\sum_{t\le k}|c^*_{i,t}|\le M_L:=k+1+\frac1{|L_0|}
\end{equation}
(\lean{exact\_relation}, \lean{abs\_eta\_le}, \lean{slope\_norm\_bound}, \lean{nLogSlopeBound}). The
correction is always made at the logarithm of index $0$; this is why $L_0\ne0$ and $k\ge1$ are needed.
The bound $M_L$ is uniform in $i$; this is the reason for the normalization $|r_{i,j}|\le1$ of the good
pairs. By \eqref{exact}, $\ell_i^*(JL)=J\varrho_i$: the centers \eqref{centers} are the points
$\big(e^{JL};\ell^*(JL)\big)$.

\subsection{The translation identity and the row expansion}
Files \lean{CLogSlopeTranslation}, \lean{CLogFullSeriesBridge}, \lean{CLogActualSlopeMatrix};
\lean{NLogSlopeTranslation}, \lean{NLogRowExpansion}, \lean{NLogSeriesFlatten},
\lean{NLogFullSeriesBridge}, \lean{NLogMainPeriodBridge}. For
$a\in\Z_{\ge0}^m$ and a column $P=Y^hX^\gamma$ define the entire function of $z\in\C^{k+1}$
\begin{equation}\label{entire}
f_{a,P}(z)=[u^a]\,P\big(e^{z_0},\dots,e^{z_k},(\ell_i^*(z)+u_i)_i\big)
=e^{h\cdot z}\prod_i\binom{\gamma_i}{a_i}\ell_i^*(z)^{\gamma_i-a_i}.
\end{equation}
Put $z_J(t)=JL+\log(1+t)$. Then $e^{z_J(t)}=(\beta_j^J(1+t_j))_j$, and by \eqref{exact}
$J\varrho_i+\ell_i(\log(1+t))=\ell_i^*(z_J(t))+\eta_i\log(1+t_0)$. The binomial formula in the additive variables
gives the finite identity
\begin{equation}\label{translation}
[u^\beta]\,P\big(b_J(1+t),(J\varrho_i+\ell_i(\log(1+t))+u_i)_i\big)
=\sum_{a\ge\beta}\binom a\beta\,\eta^{a-\beta}\,\log(1+t_0)^{|a-\beta|}\,f_{a,P}(z_J(t)),
\end{equation}
in which only $a$ with $\sum_iw_ia_i\le H$ occur. The correction depends on $t_0$ only. Extracting
$[t^\lambda]$: the entry of $A_H^{\rm full}$ in the row $(J,\lambda,\beta)$ and the column $P$ is
(\lean{slopeEntryK\_eq\_sum})
\begin{equation}\label{rowexp}
\sum_{a\ge\beta}\ \sum_{0\le\ell\le\lambda_0}s(a,\beta,\ell)\cdot[t^{\lambda-\ell e_0}]\,f_{a,P}(z_J(t)),\qquad
s(a,\beta,\ell)=\binom a\beta\,\eta^{a-\beta}\,[x^\ell]\log(1+x)^{|a-\beta|},
\end{equation}
where $e_0$ is the first unit vector. The scalars $s$ do not depend on the column. The pairs
$(a,\ell)$ are the \emph{choices} of the row (\lean{RowChoice}); their number is at most
$(\lfloor H\rfloor+1)^{2m}(\lfloor H/v\rfloor+1)$ (\lean{rowChoice\_card\_le}). An \emph{assignment}
gives every row one of its choices (\lean{Assignment}). By multilinearity
(\lean{det\_slopeEntryK\_eq\_sum}, \lean{slopeMatrix\_det\_eq\_assignment\_sum}), and by
Lemma~\ref{lem:shear} (\lean{actualMinor\_det\_eq\_slopeMatrix}),
\begin{equation}\label{assign}
\Delta_H=\sum_{\rm assignments}\ \prod_{\rm rows}s(a,\beta,\ell)\cdot
\det\Big([t^{\lambda-\ell e_0}]\,f_{a,P}(z_J(t))\Big)_{\rm rows,\ selected\ columns},
\end{equation}
with $(a,\ell)$ the choice assigned to the row.

The functions $f_{a,P}$ are defined by the forms $\ell_i^*$, which are homogeneous linear forms with
the bounds \eqref{exact}; the constants $\varrho_i$ do not occur in \eqref{entire}, \eqref{rowexp} and
\eqref{assign}. They are present only in the fact that $\ell_i^*$ does not vanish at the points $JL$.
The estimates of the rest of this section never use such a vanishing: they hold for the Taylor
coefficients of the functions \eqref{entire} at arbitrary points of a polydisc
(\lean{main\_period\_coeff\_eq\_mainRowK}, Lemma~\ref{lem:det}).

\subsection{Scalars, growth and row tests}
Let $c_L=\max(1,1/|L_0|)$ (\lean{nLogScalarConstant}). By \eqref{exact} and $\log q_i>w_i-1$,
$|\eta_i|\le c_Le^{\nu-\nu w_i}$ (\lean{norm\_eta\_le}). On $|x|=\frac12$ we have $|\log(1+x)|\le\log2<1$.
With $\binom a\beta\le2^{|a|}$ and $|a|\le H/w_*$,
\begin{equation}\label{scalar}
|s(a,\beta,\ell)|\le2^\ell\exp\Big(-\nu\sum_iw_i(a_i-\beta_i)+\frac{(\nu+\log2+\log c_L)H}{w_*}\Big).
\end{equation}
Let (\lean{nLogRadius}, \lean{nLogAdditiveBound})
\[
R=4K\Big(1+\sum_{t\le k}|L_t|\Big),\qquad B_R=2R\,M_L .
\]
Expand $f_{a,P}(z)=\sum_{p\in\Z_{\ge0}^{k+1}}c_p(a,P)\prod_t(z_t/R)^{p_t}$. Then
(\lean{norm\_columnTaylorCoeff\_le})
\begin{equation}\label{growth}
|c_p(a,P)|\le\exp\Big(H\Big(\frac R{\omega}+\frac{\log B_R}{w_*}\Big)\Big)\qquad\text{for all }p .
\end{equation}
In the formal proof this is obtained without Cauchy's inequalities: each of the $|\gamma-a|$ linear
factors in \eqref{entire} is expanded as a sum over the choice of a variable $z_t$
(\lean{NLogFormExpansion}), which writes $f_{a,P}$ as a finite sum of terms
$\mathrm{coef}\cdot\prod_tz_t^{e_t}e^{h_tz_t}$ with $\sum|\mathrm{coef}|\le2^{|\gamma|}M_L^{|\gamma-a|}$;
the normalized Taylor coefficients of such a term are at most $R^{|e|}e^{R|h|}$; and $|h|\le H/\omega$,
$|\gamma|\le H/w_*$. Further, $|JL_t|+\frac34\le R/2$ for $J<K$ (\lean{nLog\_center\_radius\_bound}), and
the \emph{row test} $[t^{\lambda'}]$ applied to the normalized monomial of index $p$ at the center $J$
satisfies (\lean{norm\_tupleTest\_le})
\begin{equation}\label{test}
\Big|[t^{\lambda'}]\prod_t\big(z_{J,t}(t)/R\big)^{p_t}\Big|\le2^{|\lambda'|}\,2^{-|p|}.
\end{equation}
With the factor $2^\ell$ of \eqref{scalar} the total cost of the row tests is $2^{|\lambda|}\le
e^{H\log2/v}$ per row.

\subsection{Collisions of Taylor coefficients in \texorpdfstring{$k+1$}{k+1} variables}
Files \lean{NLogCollisionCombinatorics}, \lean{NLogCollisionGeometricSum}, \lean{NLogFiniteCollision},
\lean{NLogCollisionExponential}, \lean{NLogCollision}, \lean{NLogAnalyticSeries},
\lean{NLogAnalyticColumns}.

\begin{lemma}[\lean{degree\_sum\_lower}]\label{lem:coll}
For $N$ distinct vectors $p\in\Z_{\ge0}^{k+1}$,
\[
\sum|p|\ \ge\ \tfrac14N^{\frac{k+2}{k+1}}-N .
\]
\end{lemma}
\begin{proof}
Let $\tau=\lfloor\frac12N^{1/(k+1)}\rfloor$. The vectors with $|p|<\tau$ lie in the cube
$\{0,\dots,\tau-1\}^{k+1}$, so there are at most $\tau^{k+1}\le N/2^{k+1}\le N/2$ of them. The other
vectors, at least $N/2$ in number, have $|p|\ge\tau>\frac12N^{1/(k+1)}-1$.
\end{proof}

The constant $\frac14$ is not optimal but it is the same for all $k$; only the exponent
$(k+2)/(k+1)$ matters.

\begin{lemma}[\lean{weighted\_main\_row\_determinant\_exp\_bound}]\label{lem:det}
Fix an assignment, and let $n_a$ be the number of rows whose choice has first component $a$. Then the
determinant in \eqref{assign} is at most
\[
\exp\Big(-c_0\sum_an_a^{\frac{k+2}{k+1}}+MH\Big(\frac R\omega+\frac{\log B_R}{w_*}+\frac{\log2}{v}+\mathrm{rem}_H\Big)\Big)
\cdot\prod_{\rm rows}2^{-\ell},\qquad c_0=\frac{\log2}8
\]
in absolute value, where $\mathrm{rem}_H=\big(\log M+\frac12\log2-(k+1)\log(1-2^{-1/2})\big)/H$
(\lean{collisionConstant}, \lean{collisionRemainder}).
\end{lemma}
\begin{proof}
Expand every $f_{a,P}$ in normalized monomials. By multilinearity the determinant is a sum over the
assignments of an index $p$ to every row; the passage to infinite sums is made through finite sets of
indices (\lean{NLogFiniteCollision}, \lean{infinite\_collision\_bound}). The row of coefficients
$(c_p(a,P))_P$ depends on $(a,p)$ only, so a term in which two rows have the same pair $(a,p)$ vanishes
(\lean{coefficient\_det\_eq\_zero\_of\_collision}). In the other terms the indices within every
$a$-block are distinct, and Lemma~\ref{lem:coll} gives $\sum_{\rm rows}|p|\ge\frac14\sum_an_a^{(k+2)/(k+1)}-M$.
Write the factor $2^{-|p|}$ of \eqref{test} as $s^{|p|}s^{|p|}$ with $s=2^{-1/2}$. The product of the
first halves is at most $s^{\frac14\sum_an_a^{(k+2)/(k+1)}-M}$; the second halves are summed over all
indices and contribute $(1-s)^{-(k+1)}$ per row. A coefficient determinant is at most $M!$ times the
$M$-th power of \eqref{growth}, and the row tests contribute $\prod2^{|\lambda|-\ell}\le
e^{MH\log2/v}\prod2^{-\ell}$ (\lean{collision\_factor\_le\_exp\_normalized}).
\end{proof}

\subsection{The dichotomy}
Files \lean{NLogGradeDichotomy}, \lean{CLogGradeData}, \lean{CLogAnalyticAssignment},
\lean{CLogAnalyticSlopeBound}, \lean{CLogAnalyticData}, \lean{CLogDataCounting}. Let
$D_A(H)$ be the number of $a\in\Z_{\ge0}^m$ with $\sum_iw_ia_i\le AH$ (\lean{lowIndexCount}); then
$D_A(H)\sim A^mH^m/(m!\prod_iw_i)$. Consider one assignment.

If at least $\eta M$ rows are assigned an $a$ with $\sum_iw_ia_i\le AH$, then by H\"older's inequality
(\lean{occupancy\_mass\_lower})
$\sum_an_a^{(k+2)/(k+1)}\ge(\eta M)^{(k+2)/(k+1)}D_A(H)^{-1/(k+1)}$. Otherwise more than $(1-\eta)M$ rows have
$\sum_iw_ia_i>AH$, and $\sum_{\rm rows}\sum_iw_i(a_i-\beta_i)\ge MH\big(A(1-\eta)-\overline b_H\big)$
(\lean{occupancy\_grade\_dichotomy}). In the first case bound the exponential in \eqref{scalar}
involving $a-\beta\ge0$ by $1$; in the second drop the collision term. With
\begin{equation}\label{rate}
\Gamma_H=\frac{c_0\,\eta^{\frac{k+2}{k+1}}}{H}\Big(\frac{M_H}{D_A(H)}\Big)^{\frac1{k+1}}
\ \longrightarrow\ \Gamma=c_0\,\eta^{\frac{k+2}{k+1}}
\Big(\frac{K\theta^m}{(m+1)(m+2)\cdots(m+k+1)\,v^{k+1}A^m}\Big)^{\frac1{k+1}}
\end{equation}
(\lean{collisionRate}, \lean{collisionLimit}, \lean{risingFactor}, \lean{tendsto\_collisionRate}; the
limit follows from \eqref{count}), every term of \eqref{assign} is at most
$\exp\big(MH(E_{\rm an}+\mathrm{rem}_H+\max\{-\Gamma_H,-\nu(A(1-\eta)-\overline b_H)\})\big)$
(\lean{assignmentTerm\_norm\_le}), where
\begin{equation}\label{Ean}
E_{\rm an}=\frac R\omega+\frac{\log2}{v}+\frac{\log B_R+\nu+\log2+\log c_L}{w_*}
\end{equation}
(\lean{analyticError}, \lean{nLogWeightErrorCoefficient}). The number of assignments is at most
$\big((\lfloor H\rfloor+1)^{2m}(\lfloor H/v\rfloor+1)\big)^{M}$. Summing:

\begin{proposition}[\lean{CLogData.actualMinor\_analytic\_upper\_bound}]\label{prop:upper}
Let $\nu\ge0$. For every integer $H>0$ and every selection with $\Delta_H\ne0$,
\[
\frac{\log|\Delta_H|}{M_HH}\ \le\ E_{\rm an}+\mathrm{rem}'_H+\max\big\{-\Gamma_H,\ -\nu\big(A(1-\eta)-\overline b_H\big)\big\},
\]
where $\mathrm{rem}'_H=\mathrm{rem}_H+\log\big((\lfloor H\rfloor+1)^{2m}(\lfloor H/v\rfloor+1)\big)/H\to0$ as
$H\to\infty$ (\lean{analyticRemainder}, \lean{tendsto\_analyticRemainder}).
\end{proposition}

\section{Parameters and the end of the proof}\label{sec:par}

Files \lean{NLogParameters}, \lean{NLogDimension}, \lean{NLogMargins}, \lean{NLogFiniteParameters},
\lean{CLogData}, \lean{CLogDataConstruction}, \lean{CLogComparison}, \lean{CLogEndpointReduction},
\lean{CLogMainTheorem}.

\subsection{The balance}

\begin{lemma}[\lean{exists\_rational\_balance}]\label{lem:balance}
Let $n\ge1$ and $\nu>n+1$. There are rational numbers $\theta,A,B,C$ with
\[
0<\theta<A<B<1<C,\qquad1-\theta<\nu(A-\theta),\qquad C^nB<1,\qquad B<C\theta<1 .
\]
\end{lemma}
\begin{proof}
Choose a rational $\mu$ with $(n+1)/\nu<\mu<1$ and put
\[
\lambda=\frac{1+\mu}2,\qquad s=1-\frac{1-\lambda}{n},\qquad c=1-(1-\lambda)(1-s),\qquad
\theta=s^{n+1},\qquad A=c\,s^{n}.
\]
Then $\mu\le\lambda^2$, $\lambda\le s<c<1$, $c-s=\lambda(1-s)$, and $s^n\ge\lambda$ by Bernoulli's
inequality. Hence $\theta<A<1$,
\[
1-\theta\le(n+1)(1-s)<\nu\mu(1-s)\le\nu\lambda s^n(1-s)=\nu(A-\theta),
\]
and $A^{n+1}=c^{n+1}\theta^n<\theta^n$. The last inequality gives $(A/\theta)^nA<1$, and
$(A/\theta)\theta=A<1$; so there is a rational $C>A/\theta>1$ with $C^nA<1$ and $C\theta<1$. Then
$A<\min(C^{-n},C\theta)$, and $B$ is a rational number between $A$ and this minimum.
\end{proof}

The inequality $A^{n+1}<\theta^n$ together with $1-\theta<\nu(A-\theta)$ is possible exactly because
$\nu>n+1$; this is where the exponent of Theorem~\ref{main} comes from.

\begin{definition}[\lean{Parameters k nu}]\label{def:par}
\emph{Parameters} for $k$ and $\nu$ are rationals $\theta,A,B,C$ as in Lemma~\ref{lem:balance} with
$n=k+1$, and a real number $\eta\in(0,1)$ with
\begin{equation}\label{gap}
g:=\nu\big(A(1-\eta)-\theta\big)-(1-\theta)>0 .
\end{equation}
\end{definition}
Parameters exist for $\nu>k+2$ (\lean{exists\_parameters}).

\subsection{The dimension}
For parameters and an integer $t\ge1$ put (\lean{dimensionK}, \lean{dimensionW}, \lean{dimensionV},
\lean{ThreeLogSigmaBounds.sigma})
\begin{equation}\label{values}
m=(k+1)t,\qquad K=\lfloor C^{(k+1)t}\rfloor,\qquad\omega=B^{-t},\qquad v=2\Big(\frac{C\theta}{B}\Big)^{t},
\qquad\sigma=\frac1{100(m+k+4)} .
\end{equation}
These are rational. We have $\omega/v=\frac12(C\theta)^{-t}$ and
\begin{equation}\label{stages}
K\Big(\frac{\omega}{v}\Big)^{k+1}\theta^{m}<\frac12,\qquad
\Big(\frac{\omega}{v}\Big)^{j}\theta^{e}<\frac12\quad(j\ge1,\ e\ge jt),\qquad
\frac K\omega\le(C^{k+1}B)^t
\end{equation}
(\lean{total\_lt\_half}, which uses $k\ge1$; \lean{stage\_lt\_half}; \lean{dimension\_K\_div\_W\_le}).

\begin{lemma}[\lean{margins\_of\_pos\_bounds}]\label{lem:margins}
Let $k\ge1$, $k+1\le m$, and suppose $(1+3\sigma)^e<2$ for $e\le m+k+1$,
\[
K\Big(\frac\omega v\Big)^{k+1}\theta^m<\frac12,\qquad
\Big(\frac\omega v\Big)^{k+1-s}\theta^{m-s}<\frac12\ \ (1\le s\le k-1),\qquad
\frac\omega v\,\theta^{m-k}<\frac12,\qquad\theta^{m-k}<\frac12,
\]
and $(k+3)(m+k+2)/\sigma<\omega/v$. Then $\mathrm{Mar}_k(m,K;\omega,v,\theta,\sigma)$ holds.
\end{lemma}
This is the unrolled form of Definition~\ref{def:margins}: the bound with index $s$ is the volume
margin after $s$ fibre reductions.

\begin{lemma}[\lean{exists\_baseData}]\label{lem:data}
Let $k\ge1$ and let parameters be given. For all real $b_1$, $b_0$, $\kappa\ge0$, $\eps'>0$ and
$\Gamma_0$ there is $t$ such that \eqref{values} are data for $k$ slopes and
\begin{equation}\label{first}
\frac{b_1m+b_0}{v}+\frac{\kappa\,K}{\omega}<\eps',\qquad
\Gamma_0<\eta^{\frac{k+2}{k+1}}\Big(\frac{K\theta^m}{(m+1)\cdots(m+k+1)\,v^{k+1}A^m}\Big)^{\frac1{k+1}} .
\end{equation}
\end{lemma}
\begin{proof}
For large $t$: $(1+\sigma)\theta<1$; the inflation factors are less than $2$ because
$\sigma=1/(100(m+k+4))$; the volume bounds of Lemma~\ref{lem:margins} follow from \eqref{stages}, with
$(j,e)=(k+1-s,m-s)$ and $(1,m-k)$, and from $\theta^{m-k}\to0$; the torus margin holds because its left
side is at most $400(m+k+4)^3$ while $\omega/v=\frac12(C\theta)^{-t}$ grows exponentially
(\lean{tendsto\_torus\_error}). The first inequality in \eqref{first} holds because $v\to\infty$ and
$K/\omega\to0$ exponentially in $t$ while $m$ is linear in $t$ (\lean{tendsto\_dimension\_error}). For the second, $v^{k+1}=2^{k+1}(C\theta/B)^m$ and $K\ge\frac12C^m$ for large
$t$, so the quantity in brackets is at least $(B/A)^{(k+1)t}/(2^{k+2}(m+k+1)^{k+1})$, and its
$(k+1)$-th root is at least $(B/A)^t/(4(m+k+1))$, which tends to infinity
(\lean{collision\_minorant\_le}, \lean{tendsto\_collision\_minorant}).
\end{proof}

\subsection{The margins of the setting}
In the setting of \ref{sec:setting} the following three inequalities are required (fields
\lean{dimension\_margin}, \lean{collision\_margin}, \lean{weight\_margin} of \lean{CLogData}):
\begin{gather}
\frac{\Lambda_0F_0+\log2}{v}+\frac{\big(4(1+\sum_t|L_t|)+\sum_t\log\operatorname{den}(\beta_t)\big)K}{\omega}<\frac{\eps_0}2,
\label{m1}\\
\Gamma>2,\qquad\qquad\frac{\log B_R+\nu+\log2+\log c_L}{w_*}<\frac{\eps_0}2 ,\label{m2}
\end{gather}
with $\Gamma$ from \eqref{rate}. They give (\lean{error\_sum\_lt\_epsilon},
\lean{collision\_exceeds\_error\_sum})
\begin{equation}\label{E}
E:=E_{\rm ar}+E_{\rm an}<\eps_0<g,\qquad\Gamma>2>1+E .
\end{equation}

\begin{proposition}[\lean{CLogComparison.actual\_minor\_bounds\_inconsistent}]\label{prop:cmp}
Let $\nu>1$. In the setting of \ref{sec:setting}, with \eqref{m1} and \eqref{m2}, it is impossible that
for every $H_1$ there are an integer $H\ge H_1$ and a selection with $\Delta_H\ne0$.
\end{proposition}
\begin{proof}
For such $H$, Propositions~\ref{prop:lower} and~\ref{prop:upper} give
\[
-(1-\overline b_H)-E_{\rm ar}\le E_{\rm an}+\mathrm{rem}'_H+\max\big\{-\Gamma_H,\ -\nu\big(A(1-\eta)-\overline b_H\big)\big\}.
\]
If the maximum is the second term, then $\nu A(1-\eta)-1-(\nu-1)\overline b_H\le E+\mathrm{rem}'_H$; the left
side decreases in $\overline b_H\le\theta$ and is at least $g$, so $g\le E+\mathrm{rem}'_H$. If it is the
first, then $\Gamma_H\le1-\overline b_H+E+\mathrm{rem}'_H\le1+E+\mathrm{rem}'_H$. Since $\mathrm{rem}'_H\to0$ and
$\Gamma_H\to\Gamma$, both contradict \eqref{E} for large $H$.
\end{proof}

\subsection{Proof of Proposition \ref{prop:no}}
Let $k\ge1$, assume $\Unif(k)$, and let $a_0,\dots,a_{k-1},g$ and $\nu>k+2$ be as in the proposition.
Suppose a $\Phi_k$-supply of $\nu$-good pairs exists. Put $\beta_j=a_j$ $(j<k)$, $\beta_k=g$; then
$L_0\ne0$ (\lean{snocLogs\_zero\_ne}; here $k\ge1$ is used) and the $\beta_t$ have no power relation
(\lean{noPowerRelation\_snoc}). The choices are made in the following order
(\lean{exists\_cLog\_data}, \lean{no\_supplyC}).
\begin{enumerate}[leftmargin=2em]
\item By $\Unif(k)$, fix thresholds $T(D)$ for all data $D$ (Lemma~\ref{lem:minors}).
\item Choose parameters $\theta,A,B,C,\eta$ for $k$ and $\nu$, and put $\eps_0=\min(g/2,\frac12)$.
\item Choose $t$ by Lemma~\ref{lem:data} with $b_1=0$, $b_0=\Lambda_0F_0+\log2$,
  $\kappa=4(1+\sum_t|L_t|)+\sum_t\log\operatorname{den}(\beta_t)$, $\eps'=\eps_0/2$ and $\Gamma_0=2/c_0$.
  This fixes the data $D$, hence $m$, $K$, $\omega$, $v$, $\sigma$, $R$, $B_R$, and gives \eqref{m1} and
  $\Gamma>2$.
\item Choose $X\ge0$ with $(\log B_R+\nu+\log2+\log c_L)/(X+1)<\eps_0/2$, and put $w_*=X+1$.
\item Apply Lemma~\ref{lem:sel} with $m$, the thresholds $T(D)$ and the bound $X+1$. This gives
  $\nu$-good pairs $(r_i,\varrho_i)$, $0\le i<m$, with slope vectors in general position, $q_i\ge2$, weights
  $w_i=\lceil\log q_i\rceil$ above $T(D)$, and $w_i>w_*$; so \eqref{m2} holds.
\item By Lemma~\ref{lem:minors} there are arbitrarily large integers $H$ with a selection and
  $\Delta_H\ne0$. This contradicts Proposition~\ref{prop:cmp}.
\end{enumerate}
There is no simultaneous limit in the dimension, the denominators and the degree. \qed

\subsection{Proof of Theorems \ref{step} and \ref{main}}
Theorem~\ref{step} follows from Propositions~\ref{prop:nosupply} and~\ref{prop:no}. Theorem~\ref{main}
holds for $n=1$ by \cite{LJZ,L} (\lean{measureC\_one}) and for $n\ge2$ by Theorems~\ref{step}
and~\ref{thm:U} with $k=n-1$, by induction on $n$. Corollary~\ref{cor:hom} is the case $h_0=0$ of
Theorem~\ref{main}, since $H(0,h)=\max_i|h_i|$ for $h\ne0$. \qed
\section{Remarks}\label{sec:rem}

\subsection{What is used from \texorpdfstring{\cite{L}}{[L]}}
The method of \cite{L} is the separated-weight interpolation determinant of \cite{OAI}: weighted
simplices of monomials and of jets, thresholds for the weights chosen before the centers, and the
comparison of an arithmetic lower bound with an analytic upper bound for one nonzero minor. From \cite{L}: the irrationality of $\log a$ and the theorem that its irrationality exponent is $2$,
for positive rational $a\ne1$ (Lemma~\ref{lem:aff} and the case $n=1$; the latter is due to
\cite{LJZ}); the two-torus curve inequality (base of the induction in Section~\ref{sec:curve}); the
general parts of the geometric argument, which do not refer to the coordinates (the transverse
multiplicity lemma for an arbitrary commuting frame, normal bases, the blow-up, nefness and
ampleness, Serre vanishing, the Rees algebra); the successive thresholds; lattice point counts in
weighted simplices; the bound $\log\operatorname{lcm}(1,\dots,n)\le\Lambda_0n$; the denominator-clearing
lemma; the determinant lemma for weighted row sums and the count of choices. Everything else is
proved here for an arbitrary number of torus coordinates.

Two results of \cite{L} are not used and are obtained again: the two-torus interpolation lemma
(Theorem~\ref{thm:U} for $k=1$ contains it, without the restriction to nonzero slopes) and the theorem
on linear forms in two logarithms (Corollary~\ref{cor:hom} for $n=2$). Appendix~\ref{app:layers} gives
the list of the formal statements on which the proof of Theorem~\ref{main} depends.

\subsection{Statements that are not formalized}
The optimality of the exponent (Minkowski) and the comparison with Baker's theory in Section~1 are
context and are not formalized. The following consequences of
Theorem~\ref{main} are standard and are not formalized either.

\emph{Arbitrary positive rational numbers.} For positive rationals $a_1,\dots,a_n$, with or without
power relations, and every $\eps>0$ there is $C>0$ such that for all $(h_0,h)\in\Z\times\Z^n$ either
$h_0+\sum_ih_i\log a_i=0$ or \eqref{eq:main} holds: express the $a_i$ as monomials in a basis of the
multiplicative group which they generate and apply Theorem~\ref{main} to the basis.

\emph{Simultaneous approximation.} By Khintchine's transference principle \cite[Chapter~V]{Ca},
Theorem~\ref{main} is equivalent, up to $\eps$, to the following statement: for positive rationals
$a_1,\dots,a_n$ without a power relation and every $\eps>0$ there is $c>0$ such that
\[
\max_{1\le i\le n}|q\log a_i-p_i|\ \ge\ c\,q^{-1/n-\eps}
\]
for all integers $q\ge1$ and $p_1,\dots,p_n$. The exponent $1/n$ is again best possible, by
Dirichlet's theorem.

\subsection{Ineffectivity}
The constant $C$ of Theorem~\ref{main} is not effective. For each $n\ge2$ the proof shows that there
cannot be $m$ good pairs with slope vectors in general position and with denominators growing as fast
as the thresholds require, where $m=nt$ depends on $\eps$; it gives no bound for the denominators of
the good pairs themselves. In addition the induction uses the constants of $\Meas(n-1)$ in
Lemma~\ref{lem:plane}, and its base is the theorem for $n=1$ of \cite{LJZ,L}.

\appendix
\section{The formalization}\label{app:lean}

\subsection{Statements}
The development is written in Lean~4 (version 4.34.1) with Mathlib, on top of the Lean supplement of
\cite{L}; it is in the ancillary files of this paper (directory \lean{anc/lean}).
The files \lean{QuadraticPi/NLogTarget.lean}, \lean{QuadraticPi/CLogTarget.lean} and
\lean{QuadraticPi/CLogMainTheorem.lean} contain (with the line breaks changed):
\begin{alltt}\small
def NoPowerRelation \{n : ℕ\} (a : Fin n → ℝ) : Prop :=
  ∀ h : Fin n → ℤ, ∏ i, a i ^ h i = 1 → h = 0

def height \{n : ℕ\} (h : Fin n → ℤ) : ℝ :=
  max 1 ((Finset.univ.sup fun i => (h i).natAbs : ℕ) : ℝ)

def heightC \{n : ℕ\} (h0 : ℤ) (h : Fin n → ℤ) : ℝ := max |(h0 : ℝ)| (height h)

def MeasureC (n : ℕ) : Prop :=
  ∀ a : Fin n → ℚ, (∀ i, 0 < (a i : ℝ)) → NoPowerRelation (fun i => (a i : ℝ)) →
    ∀ ε : ℝ, 0 < ε → ∃ C : ℝ, 0 < C ∧
      ∀ (h0 : ℤ) (h : Fin n → ℤ), (h0 ≠ 0 ∨ h ≠ 0) →
        C * heightC h0 h ^ (-(n : ℝ) - ε) ≤
          |(h0 : ℝ) + ∑ i, (h i : ℝ) * Real.log (a i : ℝ)|

theorem measureC_succ_all (d : ℕ) : MeasureC (d+1) := by
  induction d with
  | zero => exact measureC_one
  | succ d ih =>
    have h := CLogEndpointReduction.measureC_succ_of_uniformStatementC (k := d+1)
      (by omega) ih (CLogUniformInterpolation.cLogUniformInterpolation (by omega))
    exact h

theorem measureC_all (n : ℕ) (hn : 1 ≤ n) : MeasureC n

theorem linear_forms_in_logarithms_with_constant (n : ℕ) (hn : 1 ≤ n)
    (a : Fin n → ℚ) (ha : ∀ i, 0 < (a i : ℝ))
    (hind : NoPowerRelation (fun i => (a i : ℝ))) (ε : ℝ) (hε : 0 < ε) :
    ∃ C : ℝ, 0 < C ∧ ∀ (h0 : ℤ) (h : Fin n → ℤ), (h0 ≠ 0 ∨ h ≠ 0) →
      C * heightC h0 h ^ (-(n : ℝ) - ε) ≤
        |(h0 : ℝ) + ∑ i, (h i : ℝ) * Real.log (a i : ℝ)|
\end{alltt}
The index set \lean{Fin n} is $\{0,\dots,n-1\}$; \lean{heightC h0 h} is $\max(1,|h_0|,|h_1|,\dots,|h_n|)$.
The file \lean{CLogCorollaries} contains the implication from \lean{MeasureC n} to the proposition
\lean{Measure n} of \lean{NLogTarget} (the same statement without $h_0$, for $h\ne0$) and
Corollary~\ref{cor:hom}.

\subsection{Audit}
The file \lean{CLogAudit} runs \lean{\#print axioms} on \lean{measureC\_all},
\lean{linear\_forms\_in\_logarithms\_with\_constant}, the two statements of \lean{CLogCorollaries},
\lean{measureC\_one}, \lean{measureC\_succ\_of\_uniformStatementC}, \lean{no\_supplyC},
\lean{cLogUniformInterpolation}, the curve inequalities \lean{curve\_inequality\_one},
\lean{weighted\_curve\_inequality\_chain} and \lean{weighted\_curve\_inequalityC\_chain}, and five
further intermediate statements. For each of them the answer is
\begin{quote}
\lean{[propext, Classical.choice, Quot.sound]},
\end{quote}
the three standard axioms of Lean. The 68 files \lean{QuadraticPi/CLog*.lean} have 7\,928 lines; none
contains \lean{sorry}, \lean{admit}, an added axiom or \lean{native\_decide}. Theorem~\ref{main} and
Corollary~\ref{cor:hom} are proved in Lean with no hypotheses other than those of their statements;
no classical theorem enters as an assumption.

\subsection{The two layers, and what the proof uses}\label{app:layers}
The Lean project contains the supplement of \cite{L} (files \lean{TwoLog*}, \lean{OAI/*} and other files
of \lean{QuadraticPi}; the files \lean{OAI/*} are a selected part of the proof sources \cite{OAIL} with
local extensions), a first layer of files \lean{NLog*} for an arbitrary number of torus
coordinates and centers on the torus, and the second layer \lean{CLog*} of this paper. The second
layer is built in four ways.
\begin{itemize}[leftmargin=2em]
\item \emph{By translation.} The chart at a center $(b,\xi)$ is defined as the chart at $(b,0)$
  composed with $\tau_\xi$ (\lean{formalJetC}), and the local objects of Section~\ref{sec:packets} at
  $(b,\xi)$ as the images under $\tau_{-\xi}$ of those at $(b,0)$. The statements on charts, jet ideals,
  localizations, colengths and the exceptional degree, and the determinant-one comparison of
  Lemma~\ref{lem:shear}, are then deduced from the first layer; no local argument is repeated. For
  tuples, the same is done by the translation $u\mapsto u^T$ of Step~9.
\item \emph{With the additive coordinates carried along.} The global geometric part of
  Section~\ref{sec:packets} and the files on the numerical data of
  Sections~\ref{sec:ar}--\ref{sec:par} are parametrized by a record (\lean{GeometricData} with the
  field \lean{shifts}; \lean{CLogData} with the field \lean{rho}). Their proofs are those of the
  first layer for the new records.
\item \emph{With constants in the analytic layer.} The slope polynomials, the period series and the
  row expansion are stated with constants $\kappa_i$ added to the linear forms
  (\lean{slopePolynomialK}, \lean{periodSeriesK}, \lean{det\_slopeEntryK\_eq\_sum},
  \lean{main\_period\_coeff\_eq\_mainRowK}); for $\kappa=0$ they are the statements of the first layer.
\item \emph{New files} for the reduction with pairs (Section~\ref{sec:sel}), for the case of two torus
  coordinates and for the induction of Section~\ref{sec:curve}.
\end{itemize}

The file \lean{CLogDependencies} computes the set of all declarations of the project on which
\lean{measureC\_all} depends through types, definitions and proofs (declarations of Lean and of
Mathlib are not followed). It consists of 11\,159 declarations (8\,799 of them are theorems, whose
proofs were read): 613 in 65 files \lean{CLog*}, 934 in 77 of the 119 files \lean{NLog*}, 346 in 46 files
\lean{TwoLog*}, 642 in 63 other files of \lean{QuadraticPi}, and 8\,610 in 813 files \lean{OAI/*}. The
project also contains an earlier development for three logarithms (files \lean{ThreeLog*}); from
it the proof uses 14 declarations, in 6 files, none of which is specific to three logarithms: an
identity for powers, the logarithmic differential of a product, a bound for a determinant with
summable entries, three elementary lemmas on rational numbers and finite sets, and the definition
of $\sigma$ in \eqref{values} with seven lemmas on it. The following statements of the project are
\emph{not} among the 11\,159: the interpolation theorems for two and for three torus coordinates in
their original form, the curve inequality for three torus coordinates, the theorems on linear forms
in two and in three logarithms, and the earlier theorem for $n$ logarithms without a constant term.

\Needspace{10\baselineskip}
\subsection{Dictionary}
The left column is the notation of this paper.

\noindent\begin{longtable}{@{}p{0.47\textwidth}p{0.49\textwidth}@{}}
\toprule
Paper & Lean\\
\midrule
\endhead
$\Meas(n)$; $H(h_0,h)$; no power relation & \lean{CLogTarget.MeasureC n}; \lean{heightC}; \lean{NoPowerRelation}\\
Corollary~\ref{cor:hom} & \lean{NLogTarget.Measure n}; \lean{CLogCorollaries.measure\_of\_measureC}\\
$k$ slopes, $k+1$ torus coordinates, $m$ additive coordinates & \lean{k}, \lean{m}; index type \lean{Fin (m+k+1)}\\
torus index $j$; additive index $i$ & \lean{torusIdx k m j}; \lean{addIdx k m i}; definition by cases \lean{idxCases}\\
independent $(1,L,G)$; $q(r,\varrho)$; $\nu$-good pair & \lean{AffIndependent}; \lean{pairDen}; \lean{GoodPair L G nu r ρ}\\
$\Phi$-supply; $\Phi_k$; off a list & \lean{Supply denP (GoodP L G nu) (AdmP Φ)}; \lean{forbidden k}; \lean{OffHyperplanes}\\
general position & \lean{NLogInterpolation.GeneralPosition}\\
center $(b,\xi)$; chart; translation $\tau_\xi$ & \lean{fullCenterC}; \lean{formalJetC r b ξ}; \lean{shiftX ξ}\\
$\mathrm{Mar}_k(m,K;\omega,v,\theta,\sigma)$ & \lean{Margins k m K w0 v theta sigma}\\
data $D$; $D^{\rm fib}$; Lemma~\ref{lem:torus} & \lean{BaseData k}; \lean{BaseData.fibreData}; \lean{BaseData.torus\_margin}\\
$W$; $V$ & \lean{sourceWeights}, \lean{degreeWeights}; \lean{jetWeights}\\
$\mathcal P_{\le H}$; packet map & \lean{WeightedPolynomials}; \lean{packetMapC}\\
$N(r)$, $p_j(r)$, $\chi_r(b)$; distinct characters & \lean{charDen}, \lean{charExp}, \lean{characterValue}; \lean{CharacterDistinct}\\
thresholds; weights above thresholds & \lean{Thresholds}; \lean{AboveThresholds}\\
$\Unif(k)$ & \lean{UniformStatementC k}\\
$c$-separated; $F(D)$; $c_k(D)$; $c^+_k(D)$ & \lean{SeparatedProducts D c w}; \lean{frameComparisonConstant}; \lean{chainConstantOne k D}; \lean{chainConstantK k D}\\
$\deg_W(u)$; centered places; $h_p$; $u^T$ & \lean{weightedDegree}; \lean{placesC}; \lean{contactC}; \lean{shiftTuple}\\
$d_l$; $\rho_l$; projected tuple & \lean{fibreScale}; \lean{reducedSlope}; \lean{projectedCoordinates}\\
frame $D_x$; words; ideals $I_s$ & \lean{polynomialFrame}; \lean{polynomialFrameWord}; \lean{frameDerivativeIdeal}\\
comparison \eqref{mult} & \lean{NormalProductComparison}\\
the setting of \ref{sec:setting} & structure \lean{CLogData k nu Lambda c}, with \lean{Lambda} $=\Lambda_0$ and \lean{c} $=c_0$\\
$\varrho_i$; additive parts $J\varrho_i$ of the centers & \lean{rho}; \lean{xi}\\
$F_0$, $T_i$, $\Lambda_0$, $\mathrm D_H$ & \lean{F0}, \lean{truncationOrders}, \lean{lcmConstant}, \lean{baseDenominator}\\
rows, columns, $M_H$, $\overline b_H$, $\Delta_H$ & \lean{Row}, \lean{Column}, \lean{actualRowCount}, \lean{actualMean}, \lean{actualMinor}\\
$\log|\Delta_H|/(M_HH)$ & \lean{normalizedLogMinor}\\
$\eta_i$, $c^*_{i,t}$, $M_L$, $c_L$, $R$, $B_R$ & \lean{eta}, \lean{exactCoeffR}, \lean{nLogSlopeBound}, \lean{nLogScalarConstant}, \lean{nLogRadius}, \lean{nLogAdditiveBound}\\
$E_{\rm ar}$, $E_{\rm an}$, $c_0$, $\Gamma_H$, $\Gamma$, $\mathrm{rem}'_H$ & \lean{arithmeticError}, \lean{analyticError}, \lean{collisionConstant}, \lean{collisionRate}, \lean{collisionLimit}, \lean{analyticRemainder}\\
parameters $\theta,A,B,C,\eta$; gap $g$; $\eps_0$ & \lean{Parameters k nu}; \lean{gap\_pos}; \lean{epsilon}\\
values \eqref{values} & \lean{exists\_baseData} (\lean{dimensionK}, \lean{dimensionW}, \lean{dimensionV}, \lean{sigma (m+k)})\\
\bottomrule
\end{longtable}

\Needspace{10\baselineskip}
\subsection{Files}
The files of the second layer; the numbers in brackets are the numbers of files, and the prefix
\lean{CLog} is omitted. The files \lean{NLog*} of the first layer are quoted in the text.

\noindent\begin{longtable}{@{}p{0.34\textwidth}p{0.62\textwidth}@{}}
\toprule
Part of the proof & Lean files\\
\midrule
\endhead
Statements; the case $n=1$; Theorem~\ref{main}; Corollary~\ref{cor:hom}; audit; dependencies & \lean{Target}, \lean{MainTheorem}, \lean{Corollaries}, \lean{Audit}, \lean{Dependencies} (5)\\
\addlinespace
Section~\ref{sec:sel}: pairs, reduction, hyperplanes, supplies & \lean{Vectors}, \lean{Conclusion}, \lean{Hyperplane}, \lean{Selection} (4)\\
\addlinespace
Section~\ref{sec:int}: centers, charts, translation, $\Unif(k)$ & \lean{Interpolation}, \lean{Translation} (2)\\
\addlinespace
Section~\ref{sec:curve}: contacts at arbitrary centers, Steps 2, 3 and 7, translation of tuples & \lean{Contact}, \lean{Auxiliary}, \lean{CurveExcess}, \lean{ShiftDegree}, \lean{CurveInequality} (5)\\
\addlinespace
Section~\ref{sec:curve}: two torus coordinates; the induction; Proposition~\ref{prop:CI} & \lean{CurveBaseOne}, \lean{CurveChain} (2)\\
\addlinespace
Section~\ref{sec:packets}: jet ideals, localization, colength, exceptional degree, by translation & \lean{JetIdeals}, \lean{JetLocalization}, \lean{JetColength}, \lean{CompactJetIdeal}, \lean{CenterIdentification}, \lean{ExceptionalDegree} (6)\\
\addlinespace
Section~\ref{sec:packets}: blow-up, degrees on curves, ampleness & \lean{GeometricData}, \lean{BlowupGeometry}, \lean{DataExceptionalDegree}, \lean{BlowupPolarization}, \lean{CurveImage}, \lean{CurveCoordinates}, \lean{CurveModel}, \lean{CurveHyperplaneDegree}, \lean{CurveModelDegrees}, \lean{CurveDegree}, \lean{CurveDegreeData}, \lean{CurveMargin}, \lean{UniformCurveMargin}, \lean{BlowupMargin}, \lean{BlowupAmple} (15)\\
\addlinespace
Section~\ref{sec:packets}: sections, surjectivity, Theorem~\ref{thm:U} & \lean{MatrixFrame}, \lean{SectionBound}, \lean{SectionPackets}, \lean{MatrixGeometry}, \lean{JetSurjectivity}, \lean{MatrixTheorems}, \lean{UniformData}, \lean{UniformInterpolation} (8)\\
\addlinespace
Section~\ref{sec:ar}: matrix, Lemma~\ref{lem:shear}, minors, Proposition~\ref{prop:lower} & \lean{Matrix}, \lean{MatrixJet}, \lean{JetDeterminant}, \lean{Arithmetic}, \lean{ArithmeticBounds}, \lean{Data}, \lean{DataCounting}, \lean{ActualMinor}, \lean{DataMinors} (9)\\
\addlinespace
Section~\ref{sec:an}: exact forms, translation identity, row expansion, dichotomy, Proposition~\ref{prop:upper} & \lean{SlopeTranslation}, \lean{FullSeriesBridge}, \lean{AnalyticDataIndices}, \lean{SlopeMatrix}, \lean{ActualSlopeMatrix}, \lean{GradeData}, \lean{AnalyticAssignment}, \lean{AnalyticSlopeBound}, \lean{AnalyticData} (9)\\
\addlinespace
Section~\ref{sec:par}: data from a supply, comparison, Theorem~\ref{step} & \lean{DataConstruction}, \lean{Comparison}, \lean{EndpointReduction} (3)\\
\bottomrule
\end{longtable}

\subsection{Points at which the formal proof determined the text}
At the following points the text uses the device of the formal proof in place of a more usual
argument.
\begin{itemize}[leftmargin=2em]
\item The statements for centers with additive coordinates are obtained from those for centers on
  the torus by translation, of polynomials and of tuples; the curve inequality for arbitrary centers
  uses the inequality for centers on the torus with the same data and all $K$ torus parts, so that
  no comparison of constants for different numbers of centers is needed.
\item The induction on the number of torus coordinates starts from two torus coordinates. The
  formalization showed that its step is uniform in the number of torus coordinates, and that the
  number of slopes enters the rest of the geometric argument only through the curve inequality.
\item The determinant bound of Lemma~\ref{lem:det} is stated for arbitrary points of the polydisc;
  that the exact forms vanish at the centers is not needed, and the identification of the matrix
  entries with row tests holds with the constants $J\varrho_i$.
\item The integrality of the cleared matrix entries is proved for arbitrary denominators $q_i$ that
  clear the slopes and the additive parts of the centers.
\item Lemma~\ref{lem:coll} is proved by a cube cutoff, with the constant $\frac14$ for all $k$, instead
  of the extremal property of the simplex.
\item The bound \eqref{growth} is obtained from the expansion of the products of linear forms over
  choices of variables, with exact $\ell^1$-norms, instead of Cauchy's inequalities; the entire
  functions \eqref{entire} are finite sums of separable terms.
\item The tail substitution in Lemma~\ref{lem:shear} is a composition of substitutions in one torus
  variable each.
\item The comparison \eqref{mult} is used with the factor $d^d/\epsilon^d$ and the single constant
  $F(D)$, not with $2\,d!\,\epsilon^{-d}$.
\item The matrix entries are identified with coefficients of the chart by splitting exponents into
  their torus part and their additive part (\lean{NLogSeriesFlatten}); the coefficient of a product of
  series in separate variables is the product of the coefficients (\lean{NLogSeparatedCoeff}).
\item The parameters of Lemma~\ref{lem:balance} are given by explicit formulas.
\item The frame directions at a point are computed from the identity ``a derivation of $\mathcal P$ is
  $\sum_xD(x)\,\partial_x$'' (\lean{derivation\_eq\_sum\_pderiv}).
\end{itemize}

\subsection{Compiler evidence and reproduction}\label{app:repro}
The compiler is Lean~4.34.1, commit \lean{5045d0056413266e57c625dcd7c365b10e377c52}, with Mathlib at
commit \lean{d13f23b723b8a846827a245b89c10fc7d3f11612} and the other dependencies at the revisions
recorded in \lean{anc/lean/lake-manifest.json}; these are the settings of the supplement of \cite{L},
unchanged. The Lean project in \lean{anc/lean} consists of the 1604 files of that supplement
(snapshot \lean{logarithms-v2-2026-10-20-Lean136}), byte for byte, and of 300 added files: the 68
files \lean{CLog*} of this paper, 119 files \lean{NLog*} and 113 files \lean{ThreeLog*}. No file of
the supplement is modified. Sizes and SHA-256 hashes of the added files are in
\lean{anc/ADDITIONS-MANIFEST.json}; \lean{anc/tools/verify\_package.py} checks both sets of files.

On 22 November 2026 the three audit modules \lean{CLogAudit}, \lean{NLogAudit} and
\lean{ThreeLogAxiomAudit}, whose import closure contains all added files except
\lean{CLogDependencies}, and the module \lean{CLogDependencies} were built from a clean copy of these
sources: 5976 build jobs, no errors. The only axioms reported are \lean{propext},
\lean{Classical.choice} and \lean{Quot.sound}. The report of this build, with the output of the
audits and the list of dependencies quoted above, is
\lean{anc/docs/build-reports/fresh\_clone\_build\_2026-11-22.txt}. It was not a from-scratch rebuild
of every dependency: Mathlib was taken from its prebuilt cache at the pinned revision, and the
modules of the project were compiled. Readers start with
\lean{anc/REPRODUCTION.txt}; the command \lean{lake build QuadraticPi.CLogAudit} in \lean{anc/lean}
builds Theorem~\ref{main} and prints the axioms of the sixteen statements of the audit.

The statements formally verified are Theorem~\ref{main}, Corollary~\ref{cor:hom},
Theorems~\ref{step} and~\ref{thm:U}, and every intermediate statement for which a Lean name is
given. The statements of Section~\ref{sec:rem} marked as not formalized, and the comparisons with
the literature, are not part of the formal verification. The files \lean{OAI/*} carry the notices
and the licence of their source (\lean{anc/lean/SOURCE-NOTICE.txt},
\lean{anc/lean/LICENSE-OAI-APACHE-2.0.txt}); the added files are under the Apache License~2.0
(\lean{anc/LICENSE-ADDED-FILES.txt}).

\end{document}